\documentclass[11pt]{amsart}
\usepackage[T1]{fontenc}
\usepackage[utf8]{inputenc}
\usepackage{lmodern,microtype}
\usepackage[top=1in, bottom=1.25in, left=0.9in, right=0.9in]{geometry}
\usepackage{amsmath,amssymb,amsthm,mathtools}
\usepackage[dvipsnames]{xcolor}
\usepackage{fancyhdr}
\usepackage{tikz,float,needspace}
\usepackage[colorlinks=true,linkcolor=MidnightBlue,citecolor=MidnightBlue,urlcolor=MidnightBlue]{hyperref}
\usepackage[foot]{amsaddr}

\numberwithin{equation}{section}
\newtheorem{theorem}{Theorem}[section]
\newtheorem{proposition}[theorem]{Proposition}
\newtheorem{lemma}[theorem]{Lemma}
\newtheorem{corollary}[theorem]{Corollary}

\newtheoremstyle{remarkitalic}
  {.5\topsep}{.5\topsep} 
  {\itshape}             
  {}                     
  {\itshape}             
  {.}{.5em}{}
\theoremstyle{remarkitalic}
\newtheorem{remark}[theorem]{Remark}

\newcommand{\R}{\mathbb R}

\DeclareMathOperator{\dist}{dist}
\DeclareMathOperator{\PV}{P.V.}
\title[Pohozaev identity for the fractional $p\,$-Laplacian ]{Boundary concentration and Pohozaev identity\\for the fractional $p\,$-Laplacian}
\author{Alejandro Ortega}
\address[A. Ortega]{Dpto. de Matem\'aticas Fundamentales, Facultad de Ciencias, UNED, 28040 Madrid, Spain}
\email{\tt alejandro.ortega@mat.uned.es}


\begin{document}
\begin{abstract}
We prove a Pohozaev-like identity for bounded weak solutions of the
fractional\break $p\,$-Laplace equation with a locally Lipschitz
nonlinearity on bounded $C^2$ domains, provided that
$s>\max\{\frac{1}{2}-\frac{1}{p},0\}$. We compare the weak equations for the solution
and its spatial deformation to show that the first-order contribution
of the energy remainder is concentrated at the boundary. The boundary
term involves the weighted trace, with a coefficient expressed by an
integral of the one-dimensional boundary profile. 
\end{abstract}
\maketitle

\section{Introduction and main results}\label{sec:introduction}

Pohozaev-like identities play a main role
in the analysis of elliptic PDEs as they provide, among others,
nonexistence and uniqueness results, monotonicity formulas or radial
symmetry of solutions. 

Starting from the seminal work of Pohozaev (cf.\cite{P65}) on the classical Laplacian case, the Pohozaev identity has been extended along the years to deal with more general operators. Guedda and V\'eron (cf. \cite{GV89}) obtained the Pohozaev identity for the $p\,$-Laplacian. The variational framework developed by Pucci and
Serrin (cf. \cite{PS86}) yields the corresponding identities for higher-order equations with Dirichlet boundary conditions.
The biharmonic equation with Navier boundary was addressed by
Van der Vorst (cf. \cite{VdV91}). A systematic treatment of polyharmonic
boundary value problems can be found in the monograph by Gazzola,
Grunau and Sweers (cf. \cite{GGS10}).\newline 
Pohozaev identities express a balance between the
energy of a solution, the contribution of the nonlinearity and a boundary
integral. For the Laplacian and the $p\,$-Laplacian with homogeneous
Dirichlet conditions, the boundary integral involves the $p$-th power of the normal
derivative $\frac{\partial u}{\partial\nu}$\break (cf. \cite{P65,GV89}). For polyharmonic equations,
it involves traces of higher-order derivatives of $u$, whose form depends
on the boundary conditions (cf. \cite[Section~7.4]{GGS10}).

In the nonlocal ($0<s<1$) linear case ($p=2$), the Pohozaev identity was
obtained in the seminal work of Ros-Oton and Serra (cf.~\cite{RS14}).
In this nonlocal case, the boundary term is expressed in terms of the
weighted trace $u/d^s$, where $d$ is the distance to the boundary
$\partial\Omega$. Djitte, Fall and Weth (cf.~\cite{DFW}) extended the
identity to general Lipschitz vector fields on bounded $C^{1,1}$ domains.
In this regard, the present work represents a nonlinear extension of the
bounded-domain identity of \cite{DFW} and also connects with the dilation argument 
by Ros-Otón and Serra (cf.~\cite{RS14}). Results covering anisotropic integro-differential
operators can be found in \cite{RSV17}.
In the nonlinear case $p\ne2$, Iannizzotto, Liu, Perera and Squassina (cf.~\cite{ILPS16})
proposed a Pohozaev identity for the fractional $p\,$-Laplacian on bounded domains, with an undetermined boundary exponent and coefficient. We establish this identity for bounded weak solutions on bounded $C^2$ domains when $s>\max\{0,\tfrac12-\tfrac1p\}$ and obtain its extension to general smooth deformation fields. Ros-Oton and Serra
(cf.~\cite{RSnonexistence}) also obtained variational inequalities and
nonexistence results under suitable regularity assumptions. On the whole
$\R^N$ and assuming Berestycki-Lions growth assumptions for the
nonlinearity $f$, Ambrosio constructs (cf.~\cite[Theorem 1.1 and 1.2]{A24})
a weak solution satisfying the identity and proves it for
$W^{s,p}(\R^N)\cap C^{1,1}(\R^N)$ solutions. These results assume
$sp<N$ and the latter also requires $s<2(p-1)/p$ when $1<p<2$. 

Our approach is based on comparing the weak equation for a solution with the equation
satisfied by its spatial deformation. This comparison allows us to
show that the first-order contribution of the energy remainder is
concentrated at the boundary and to identify the coefficient of the
resulting boundary term. The estimates leading to this limit justify
the Pohozaev identity for bounded weak solutions whenever
$s>\max\{p-2,0\}/(2p)$. The proof relies on the continuity of $u/d^s$ (cf. \cite{IMS20,IM24}), 
the interior H\"older estimates obtained in \cite{BLS18,GL24},
and the higher differentiability estimates from \cite{DKLN25}.

Let $N\ge2$, $1<p<\infty$, $0<s<1$ and let $\Omega\subset\R^N$ be a
bounded open set with $C^2$ boundary. Its outward unit normal is denoted
by $\nu$ and $d(x)=\dist(x,\partial\Omega)$ for $x\in\Omega$. The
fractional Sobolev space is defined as
\[
 W^{s,p}(\R^N)=\{u\in L^p(\R^N):[u]_{s,p}<\infty\},
\]
where the Gagliardo seminorm is given by
\[
 [u]_{s,p}^p=\iint_{\R^{2N}}
 \frac{|u(x)-u(y)|^p}{|x-y|^{N+sp}}\,dx\,dy.
\]
The space $W^{s,p}(\R^N)$ is endowed with the norm
\[
 \|u\|_{W^{s,p}(\R^N)}^p
 =\|u\|_{L^p(\R^N)}^p+[u]_{s,p}^p
\]
To take into account the homogeneous Dirichlet boundary condition, we consider the subspace
\[
 W^{s,p}_0(\Omega)=\{u\in W^{s,p}(\R^N):u=0\text{ a.e. in }\Omega^c\}.
\]
The fractional $p\,$-Laplacian is the nonlocal nonlinear operator defined as
\begin{equation}\label{operator}
 (-\Delta)^s_pu(x)=C_{N,s,p}\PV\int_{\R^N}
 \frac{|u(x)-u(y)|^{p-2}(u(x)-u(y))}{|x-y|^{N+sp}}\,dy,
\end{equation}
with $C_{N,s,p}$ given in \cite[Lemma 5.1]{DGV}, namely
\vspace{-0.5cm}

\begin{equation}\label{normalization}
 C_{N,s,p}=\frac{sp(1-s)2^{2s-2}\Gamma\left(\frac{N+sp}{2}\right)}
 {\pi^{\frac{N-1}{2}}\Gamma\left(\frac{p+1}{2}\right)\Gamma(2-s)}.
\end{equation}
When $p=2$, this constant reduces to $C_{N,s,2}=\frac{4^s\Gamma(N/2+s)}{\pi^{N/2}|\Gamma(-s)|}$ and \eqref{operator} has Fourier multiplier $|\xi|^{2s}$.
The main problem we deal with is the Dirichlet problem
\begin{equation*}\tag{$P_D$}\label{eq:dirichlet-problem}
 \begin{cases}
 (-\Delta)^s_pu=f(u) &\quad \text{in }\Omega,\\
 \mkern+51mu u=0 &\quad \text{in }\Omega^c.
 \end{cases}
\end{equation*}
We assume that the nonlinearity $f$ belongs to
$C^{0,1}_{\mathrm{loc}}(\R)$ and we denote by $F(t)=\int_0^t f(\tau)\,d\tau$.

A bounded weak solution of \eqref{eq:dirichlet-problem} is a function
$u\in W^{s,p}_0(\Omega)\cap L^\infty(\Omega)$ such that
\begin{equation}\label{weak}
 \frac{C_{N,s,p}}2
 \iint_{\R^{2N}}
 \frac{|u(x)-u(y)|^{p-2}(u(x)-u(y))(\varphi(x)-\varphi(y))}
 {|x-y|^{N+sp}}\,dx\,dy
 =\int_\Omega f(u)\varphi\,dx
\end{equation}
for all $\varphi\in W^{s,p}_0(\Omega)$. The equation \eqref{weak} above is the Euler--Lagrange equation of the functional
\[
 \mathcal F_{s,p}[u]=\frac{C_{N,s,p}}{2p}[u]_{s,p}^p-\int_\Omega F(u)\,dx.
\]
\noindent\textbf{Sketch of the proof.}
The main idea is to compare $u$ with $u_h$, a small spatial deformation of itself and compute the resulting change in energy in two ways. A change of integration variables gives the volume terms of the Pohozaev identity, while the weak equation expresses the same change as a remainder whose first-order contribution concentrates at the boundary. 

For a smooth function, a spatial variation is usually described with
$X\cdot\nabla u$, where $X$ is a vector field describing a spatial
deformation. On the contrary, we do not assume that $\nabla u$ exists. Instead, we move
points a small distance $h$, we evaluate $u$ at the moved points and
then we subtract $u$. More
precisely, let $\Phi_h(x)$ be the position at time $h$ of the trajectory
starting at $x$ and moving with velocity $X$ and set
\begin{equation}\label{deformation}
 u_h=u\circ\Phi_h,\qquad D_h=u_h-u.
\end{equation}
For example, $X(x)=x$ gives $\Phi_h(x)=e^h x$, the usual dilation. We initially choose $X$ to point strictly outward on the boundary. The support moves by the inverse flow: $u_h(x)$ can be nonzero only if $\Phi_h(x)\in\operatorname{supp}u$, or equivalently $x\in\Phi_{-h}(\operatorname{supp}u)$. Thus, for small $h>0$,
\[
 \operatorname{supp}u_h\subset\Phi_{-h}(\overline\Omega)\subset\Omega.
\]
Although $D_h$ need not vanish in a whole boundary neighborhood, because it also contains $-u$; the key point is that both terms belong to $W^{s,p}_0(\Omega)$. Hence $D_h$ is an admissible weak test. 

On the other hand, the deformed function $u_h$ satisfies a weak equation involving a
perturbed kernel. Comparison with the equation for $u$ will control
the energy remainder in the interior. To establish this control we
proceed as follows.

The function $|t|^p/p$ differs from its tangent-line approximation at $a$ by the nonnegative remainder
\begin{equation}\label{Bdef}
 B_p(a,b)=\frac{|b|^p-|a|^p}{p}-|a|^{p-2}a(b-a).
\end{equation}
Apply this expansion to the energy increments of $u$ and $u_h$ and
apply the corresponding Taylor expansion to term involving the nonlinearity
\begin{align}
 R_h&=\frac{C_{N,s,p}}2
 \iint_{\R^{2N}}
 \frac{B_p(u(x)-u(y),\,u_h(x)-u_h(y))}
 {|x-y|^{N+sp}}\,dx\,dy,\label{Rdef}\\
 T_h&=\int_\Omega\bigl(F(u_h)-F(u)-f(u)D_h\bigr)\,dx.\label{Tdef}
\end{align}

\begin{remark}
Since $B_2(a,b)=|b-a|^2/2$, in the quadratic case 
\begin{equation}\label{quadratic-remainder}
 R_h=\frac{C_{N,s,2}}4[D_h]_{s,2}^2=\frac{1}{2}\|(-\Delta)^{\frac{s}{2}}(u_h-u)\|_{L^2(\mathbb{R}^N)}^2.
\end{equation}
In the particular case $X(x)=x$, namely, the dilation $u_h(x)=u(e^h x)$, the remainder $R_h$
can be expressed in terms of the integral used by Ros-Oton and
Serra~\cite{RS14}. Set $w=(-\Delta)^{s/2}u$ and
\[
 I_\lambda=\int_{\R^N}w(\lambda x)w(x/\lambda)\,dx.
\]
Since $(-\Delta)^{s/2}u_h(x)=e^{sh}w(e^h x)$, expanding the square and the change of variable $y=e^{h/2}x$ give
\[\begin{split}
R_h
&=\frac{e^{(2s-N)h}+1}{2}\|w\|_2^2-e^{(s-N/2)h}I_{e^{h/2}}\\
&=e^{(s-N/2)h}\bigl(I_1-I_{e^{h/2}}\bigr)+\frac12\bigl(e^{(s-N/2)h}-1\bigr)^2I_1,
\end{split}\]
since $I_1=\|w\|_2^2$. The last term is $O(h^2)$. Since $(e^{h/2}-1)/h\to 1/2$,
we get
\[
 \lim_{h\downarrow0}\frac{R_h}{h}
 =-\frac12 \left.\frac{d}{d\lambda}\right|_{\lambda=1^+}I_\lambda.
\]
\end{remark}
By definition of $B_p$, we have
\[
\begin{aligned}
 \frac{|u_h(x)-u_h(y)|^p-|u(x)-u(y)|^p}{p}&=|u(x)-u(y)|^{p-2}(u(x)-u(y))
       [D_h(x)-D_h(y)]\\
 &\quad+B_p(u(x)-u(y),u_h(x)-u_h(y)).
\end{aligned}
\]
Multiply this equality by $C_{N,s,p}/(2|x-y|^{N+sp})$ and integrate over $\R^{2N}$. By \eqref{weak} with $\varphi=D_h$, the linear (in $D_h$) term equals $\int_\Omega f(u)D_h\,dx$. 
Therefore,

\[
 R_h=\frac{C_{N,s,p}}{2p}\bigl([u_h]_{s,p}^p-[u]_{s,p}^p\bigr)
 -\int_\Omega f(u)D_h\,dx,
\]
whereas the definition of $T_h$ gives
\[
 T_h=\int_\Omega(F(u_h)-F(u))\,dx-\int_\Omega f(u)D_h\,dx.
\]
Subtracting cancels the same linear term on both sides and produces
\begin{equation}\label{flow-exact}
 R_h-T_h
 =\frac{C_{N,s,p}}{2p}\bigl([u_h]_{s,p}^p-[u]_{s,p}^p\bigr)
 -\int_\Omega\bigl(F(u_h)-F(u)\bigr)\,dx.
\end{equation}
The right side is the difference of energies
$\mathcal F_{s,p}[u_h]-\mathcal F_{s,p}[u]$, while the left side contains
the remainders after using the equation.

The right side can be differentiated by moving the deformation from $u$ into the integration variables. For example, take a ball centered at zero and $X(x)=x$. Taking $\xi=e^h x$ and $\eta=e^h y$,
\[
 [u_h]_{s,p}^p=e^{(sp-N)h}[u]_{s,p}^p\qquad\text{and}\qquad \int_\Omega F(u_h)\,dx=e^{-Nh}\int_\Omega F(u)\,dx.
\]
Hence the right side of \eqref{flow-exact}, divided by $h$, converges to
\[
 \frac{C_{N,s,p}(sp-N)}{2p}[u]_{s,p}^p
 +N\int_\Omega F(u)\,dx.
\]
This explains the factors $N$ and $N-sp$ in the Pohozaev identity. For a general flow, Jacobians replace the exponential factors and give the kernel $L_X$ in \eqref{LX} below.

We have next to identify the limit of the left side of \eqref{flow-exact}. If $L$ is a Lipschitz constant of $f$ on a bounded interval containing the range of the zero extension of $u$, Taylor's formula gives
\[
 |F(b)-F(a)-f(a)(b-a)|\le\frac L2|b-a|^2,
 \qquad |T_h|\le\frac L2\|D_h\|_2^2.
\]
The increment estimates below (see Lemma \ref{lem:increments}) shows $T_h=o(h)$. Thus, the boundary contribution comes from $R_h/h$.
Near the boundary, the profile of $u$ is determined by the weighted trace
\[
 \psi=\frac{u}{d^s}\in C(\overline\Omega).
\]
At each $z\in\partial\Omega$, the ratio $u(x)/d(x)^s$ tends to $\psi(z)$ as $x$ approaches $z$ from inside $\Omega$. A displacement of size $h$ is small relative to the distance to the boundary where $d\gg h$. Comparison of the equations shows that the energy remainder there is $o(h)$. At distance $hr$ from the boundary, with $r$ fixed, the outward deformation changes the signed distance to approximately $h(r-1)$ when the normal velocity is one. The corresponding leading values of $u$ and $u_h$ are $h^s\psi(z)r_+^s$ and $h^s\psi(z)(r-1)_+^s$ respectively. Their difference can therefore contribute to the limiting profile after division by $h^s$.

The size of this contribution can already be read from \eqref{Rdef}. Take $x$ in a boundary strip of thickness $h$ and $|y-x|$ comparable to $h$. The function values have size $h^s$, so the numerator scales as $h^{sp}$ and the kernel as $h^{-N-sp}$. Per unit boundary area, the normal coordinate of $x$ contributes a factor $h$ and the separation vector $y-x$ contributes $h^N$. Their product is
\[
 h^{sp}\,h^{-N-sp}\,h\,h^N=h.
\]
Thus this region can leave a nonzero limit after division by $h$. This computation concerns only pairs with $d(x)$ and $|x-y|$ comparable to $h$; it does not yet estimate the entire double integral. To obtain an equality, we must also show that pairs deep inside $\Omega$, pairs with $|x-y|\gg h$ and pairs with $|x-y|\ll h$ do not leave another contribution.

To identify the local profile in the limit of $R_h/h$, fix $z\in\partial\Omega$ and first take unit outward normal velocity at $z$. The point $x=z-hr\nu(z)$ has signed inward distance $hr$, while $\Phi_h(x)$ has signed inward distance $h(r-1)+o(h)$. The representation $u=\psi d^s$ inside $\Omega$ and the extension of $u$ by zero outside $\Omega$ give the profile $g(t)=t_+^s$, where $t_+=\max\{t,0\}$. In these coordinates,
\[
 h^{-s}u(x)\longrightarrow\psi(z)g(r),\qquad
 h^{-s}u_h(x)\longrightarrow\psi(z)g(r-1).
\]
For a second point with rescaled normal coordinate $t$, the same calculation gives $g(t)$ and $g(t-1)$. Substitution in $B_p$ extracts the factor $|\psi(z)|^p$. After integration in the directions parallel to the tangent plane, the remaining profile integral is
\begin{equation}\label{Kdef}
 K_{s,p}=\iint_{\R^2}
 \frac{B_p(g(r)-g(t),\,g(r-1)-g(t-1))}
 {|r-t|^{1+sp}}\,dr\,dt.
\end{equation}
The constant $K_{s,p}$ depends only on $s,p$ and measures the energy remainder of a unit translation of the boundary profile. For normal velocity $X(z)\cdot\nu(z)=\ell>0$, rescaling $r,t$ shows
that the integral becomes $\ell K_{s,p}$. Multiplying the factor from integration parallel to the tangent plane by the normalization in \eqref{normalization} gives
\[
 C_{N,s,p}\frac{\pi^{\frac{N-1}{2}}\Gamma\left(\frac{1+sp}{2}\right)}{\Gamma\left(\frac{N+sp}{2}\right)}
 = C_{1,s,p}.
\]
The resulting boundary coefficient is independent of $N$ and is given by
\begin{equation}\label{cdef}
 c_{s,p}=\frac{C_{1,s,p}}{2}K_{s,p}
\end{equation}
Proposition~\ref{prop:constant} proves that this coefficient is finite and positive. Combining these local calculations gives the following candidate for the full limit,
\[
 \lim\limits_{h}\frac{R_h}{h}=
 c_{s,p}\int_{\partial\Omega}|\psi|^pX\cdot\nu\,d\sigma.
\]
Combining this limit together with $T_h=o(h)$ and the change of variables in \eqref{flow-exact}, will produce the Pohozaev identity. The main task is then to justify the limit of the integral over all pairs. On bounded rescaled regions separated from the diagonal, continuity of $\psi$ gives convergence to the profile. Nonnegativity of $B_p$ gives the lower bound for the full integral. For the upper bound, H\"older estimates control pairs near the diagonal and comparison of the equations, followed by summation over distance layers, controls pairs outside the bounded rescaled regions.

Although the proof starts with strictly outward fields, for which the deformed tests are admissible, the Pohozaev identity obtained  is linear in $X$. Writing an admissible field as the difference of two strictly outward fields and subtracting the identities for those two fields gives the general case.

For $X\in C^2(\R^N;\R^N)$ with bounded first and second derivatives, write
\begin{equation}\label{LX}
 L_X(x,y)=\operatorname{div}X(x)+\operatorname{div}X(y)
 -(N+sp)\frac{(X(x)-X(y))\cdot(x-y)}{|x-y|^2}.
\end{equation}
We can now state our main result.

\begin{theorem}\label{thm:general}
Let $N\ge2$, $1<p<\infty$ and
\begin{equation}\label{new-range}
 \max\left\{\frac12-\frac1p,0\right\}<s<1.
\end{equation}
Let $\Omega\subset\R^N$ be a bounded open set with $C^2$ boundary, let $f\in C^{0,1}_{\mathrm{loc}}(\R)$. If $u\in W^{s,p}_0(\Omega)\cap L^\infty(\Omega)$ is a weak solution of \eqref{eq:dirichlet-problem} in the sense of \eqref{weak}, its weighted trace $\psi=u/d^s$ belongs to $C(\overline\Omega)$. Moreover, for every $X\in C^2(\R^N;\R^N)$ with bounded first and second derivatives, one has
\begin{equation}\label{general-identity}
 \int_\Omega F(u)\operatorname{div}X\,dx
 -\frac{C_{N,s,p}}{2p}\iint_{\R^{2N}}
 \frac{|u(x)-u(y)|^p}{|x-y|^{N+sp}}L_X(x,y)\,dx\,dy
 =c_{s,p}\int_{\partial\Omega}|\psi|^pX\cdot\nu\,d\sigma.
\end{equation}
In particular,
\begin{equation}\label{pohozaev}
 N\int_\Omega F(u)\,dx-\frac{N-sp}{p}\int_\Omega uf(u)\,dx
 =c_{s,p}\int_{\partial\Omega}|\psi|^p x\cdot\nu\,d\sigma.
\end{equation}
\end{theorem}
For $p=2$, the coefficient agrees with that of \cite{RS14} in the normalization
\eqref{operator} and is given by
\begin{equation}\label{quadratic}
 c_{s,2}=\frac{\Gamma(1+s)^2}{2}.
\end{equation}

\textbf{Organization of the paper.}
In Section~\ref{sec:regularity} we prove the scalar estimates and establish
the regularity needed in the proof. In Section~\ref{sec:comparison} we use
that regularity to compare $u$ and $u_h$ in interior balls.
In Section~\ref{sec:boundary} we compute the boundary contribution and
complete the proof of the identity.

\section{The boundary coefficient and interior regularity}\label{sec:regularity}

We fix $q=\min\{p,2\}$ and $J(t)=|t|^{p-2}t$. Constants denoted by $C$ may vary line to line.

\begin{samepage}
\begin{lemma}\label{lem:convex}
For $a,b\in\R$,
\begin{align}
 0\le B_p(a,b)&\le C(|a|+|b|)^{p-q}|b-a|^q
 \le C(|a|^p+|b|^p),\label{Bupper}\\
 B_p(a,b)+B_p(b,a)&=(J(b)-J(a))(b-a),\label{Sidentity}\\
 |J(b)-J(a)|&\le C(|a|+|b|)^{p-q}|b-a|^{q-1},\label{Munified}\\
 (J(b)-J(a))(b-a)&\ge c(|a|+|b|)^{p-2}|b-a|^2.\label{weighted-coercivity}
\end{align}
When $1<p<2$, also
\begin{equation}\label{Msmall}
 |J(b)-J(a)|^{p/(p-1)}
 \le C(J(b)-J(a))(b-a).
\end{equation}
\end{lemma}
\end{samepage}
\begin{proof}
The function $|t|^p/p$ is convex and has derivative $J(t)$.
Its graph lies above its tangent at $a$, which proves $B_p(a,b)\ge0$.
The fundamental theorem of calculus gives
\begin{equation}\label{B-integral}
 B_p(a,b)=(b-a)\int_0^1
 [J(a+\theta(b-a))-J(a)]\,d\theta.
\end{equation}

For $p\ge2$, $J'(t)=(p-1)|t|^{p-2}$. Every point of the segment from $a$ to $b$ has absolute value at most $|a|+|b|$.
Consequently,
\[
 |J(a+\theta(b-a))-J(a)|
 \le(p-1)(|a|+|b|)^{p-2}\theta|b-a|.
\]
Integration proves the first upper bound in \eqref{Bupper} with $q=2$; taking $\theta=1$ proves \eqref{Munified}.
By symmetry, assume $|a|\ge|b|$. For $0\le\theta\le1/4$, the triangle inequality gives $|a+\theta(b-a)|\ge |a|-\theta(|a|+|b|)
 \ge\frac{|a|+|b|}{4}.$
Restricting the integral to $[0,1/4]$ gives
\[
 \begin{aligned}
 (J(b)-J(a))(b-a)
 &=(p-1)|b-a|^2\int_0^1|a+\theta(b-a)|^{p-2}\,d\theta\ge\frac{p-1}{4}|b-a|^2\left(\frac{|a|+|b|}{4}\right)^{p-2}.
 \end{aligned}
\]

For $1<p<2$, concavity of $t^{p-1}$ on $[0,\infty)$ gives
$|J(b)-J(a)|\le C|b-a|^{p-1}$ for endpoints of the same sign.
For opposite signs, use
$|a|^{p-1}+|b|^{p-1}\le2^{2-p}(|a|+|b|)^{p-1}$.
Thus this bound holds for arbitrary signs. Substitution in
\eqref{B-integral} proves \eqref{Bupper} with $q=p$ and the same bound is \eqref{Munified}. The integral formula for
$(J(b)-J(a))(b-a)$ above remains valid, the possible singularity at zero is integrable since $p-2>-1$. As $p-2<0$ and
$|a+\theta(b-a)|\le|a|+|b|$, that integral is at least
$(|a|+|b|)^{p-2}$. This proves \eqref{weighted-coercivity}.
Finally, $|J(b)-J(a)|^{1/(p-1)}\le C|b-a|$ and $J(b)-J(a)$ has the sign of $b-a$. Multiplication by
$|J(b)-J(a)|$ gives \eqref{Msmall}.

The bound by $C(|a|^p+|b|^p)$ follows from $|b-a|\le|a|+|b|$ while \eqref{Sidentity} is immediate from \eqref{Bdef}. 
\end{proof}

The next proposition establishes that the boundary coefficient in the main theorem is well defined.

\begin{samepage}
\begin{proposition}\label{prop:constant}
The integral in \eqref{Kdef} is finite and
\begin{equation}\label{strictconstant}
 K_{s,p}>\frac{2(p-1)}{sp^2}>0.
\end{equation}
For every real displacement $\ell>0$,
\begin{equation}\label{cellscale}
 \iint_{\R^2}
 \frac{B_p(g(r)-g(t),\,g(r-\ell)-g(t-\ell))}
 {|r-t|^{1+sp}}\,dr\,dt=\ell K_{s,p}.
\end{equation}
For $p=2$, \eqref{quadratic} holds.
\end{proposition}
\end{samepage}
\begin{proof}
We start by establishing the finiteness of $K_{s,p}$. It is straightforward to check that, for every $L>0$ and $\ell\in\R$, we have $g,\ g(\cdot-\ell)\in W^{s,p}((-L,L))$. Thus, by \eqref{Bupper}, we get
\[
\iint_{(-L,L)^2}\!\!\! \frac{B_p(g(r)-g(t),g(r-1)-g(t-1))}{|r-t|^{1+sp}}\,dr\,dt\le C\Bigl([g]_{W^{s,p}((-L,L))}^p+
 [g(\cdot-1)]_{W^{s,p}((-L,L))}^p\Bigr)<\infty.
\]
We continue next with the integrability at infinity. By symmetry in $r,t$ we may restrict to $r>t$. The integrand vanishes if $r\le0$ and is locally integrable on the bounded strip $0<t\le2$, $0<r-t<1$. We treat next the three remaining regions.

On $r>0\ge t$, we have $g(t)=g(t-1)=0$, so the two arguments
of $B_p$ are $r^s$ and $(r-1)_+^s$.
For $0<r\le2$, we have $0\le (r-1)_+^s\le r^s$. In particular,
the sum of the arguments is at most $2r^s$ and their difference
is at most $r^s$. Hence \eqref{Bupper} gives
$B_p(r^s,(r-1)_+^s)\le C(2r^s)^{p-q}(r^s)^q\le Cr^{sp}$.
For $r>2$, the two arguments are still bounded by $r^s$ and
their difference satisfies
\[
 r^s-(r-1)^s
 =s\int_{r-1}^r x^{s-1}\,dx
 \le s(r-1)^{s-1}
 \le s\,2^{1-s}r^{s-1}.
\]
In the last inequality we used $r-1\ge r/2$ and $s-1<0$.
Substituting these bounds into \eqref{Bupper}, we obtain
$B_p(r^s,(r-1)^s)\le C(2r^s)^{p-q}(r^{s-1})^q\le Cr^{sp-q}$.
Taking into account that $q>1$, the contribution of the region $r>0\ge t$  satisfies
\begin{align*}
 \int_0^\infty\int_{-\infty}^0
 \frac{B_p(g(r)-g(t),\,g(r-1)-g(t-1))}
      {(r-t)^{1+sp}}\,dt\,dr
 &=\frac{1}{sp}\int_0^\infty
 B_p(r^s,(r-1)_+^s)\,r^{-sp}\,dr\\
 &\quad\le C\int_0^2 r^{sp}r^{-sp}\,dr
       +C\int_2^\infty r^{sp-q}r^{-sp}\,dr\\
 &\quad=C\left(2+\int_2^\infty r^{-q}\,dr\right)<\infty.
\end{align*}
On $t>2$, $0<\tau=r-t<t/2$, the arguments of $B_p$ are
$(t+\tau)^s-t^s$ and $(t+\tau-1)^s-(t-1)^s$. Each is bounded by
$C\tau t^{s-1}$ by integrating the derivative of $t^s$ over the
corresponding interval. Similarly, the difference between the second and first
arguments is
\[
\bigl|(t+\tau-1)^s-(t-1)^s-(t+\tau)^s+t^s\bigr|
 \le C\tau t^{s-2}.
\]
Inserting these estimates into \eqref{Bupper} and using the substitution $\tau=r-t$, we obtain
\[
\begin{aligned}
 \int_2^\infty\int_t^{3t/2}
 \frac{B_p(g(r)-g(t),\,g(r-1)-g(t-1))}{(r-t)^{1+sp}}
 \,dr\,dt &\le C\int_2^\infty\int_0^{t/2}
 \frac{(\tau t^{s-1})^{p-q}(\tau t^{s-2})^q}{\tau^{1+sp}}
 \,d\tau\,dt\\
 &\quad=\frac{C}{(p-sp)2^{p-sp}}\int_2^\infty t^{-q}\,dt
 <\infty.
\end{aligned}
\]
Here $p-sp>0$ ensures integrability at $\tau=0$ and $q>1$
ensures integrability as $t\to\infty$.

On $t>0$, $\tau\ge\max\{1,t/2\}$, the increments are bounded by $C\tau^s$ by $s$-H\"older continuity. Their difference is bounded by $C(1+t)^{s-1}$ because it is a difference of unit increments at $t$ and $t+\tau$. Since $\max\{1,t/2\}\ge(1+t)/3$, using $\tau=r-t$ and applying \eqref{Bupper}, we get
\begin{align*}
\int_0^\infty\mkern-13mu\int\limits_{t+\max\{1,t/2\}}^\infty
 \mkern-13mu\frac{B_p(g(r)-g(t),\,g(r-1)-g(t-1))}
      {(r-t)^{1+sp}}\,dr\,dt&\le C\int_0^\infty\mkern-10mu\int\limits_{\max\{1,t/2\}}^\infty
 \frac{(\tau^s)^{p-q}(1+t)^{q(s-1)}}{\tau^{1+sp}}
 \,d\tau\,dt\\
 &\quad\le C\int_0^\infty(1+t)^{-q}\,dt<\infty.
\end{align*}
The bounds are finite since $q>1$. These three regions and the bounded strip $0<t\le2$, $0<r-t<1$, cover the nonzero part of the integrand for $r>t$. Symmetry in $r,t$ gives finiteness of $K_{s,p}$.

Let us prove now the lower bound \eqref{strictconstant}. On $0<r<1$ and $t<0$, the second argument of $B_p$ is zero.
The first argument is $r^s$, so $B_p(r^s,0)=(p-1)r^{sp}/p$. This region and its transpose contribute
\[
2\int_0^1\int_{-\infty}^0\frac{B_p(r^s,0)}{(r-t)^{1+sp}}\,dt\,dr=\frac{2(p-1)}p\int_0^1r^{sp} \left[\frac{(r-t)^{-sp}}{sp}\right]_{t=-\infty}^{t=0}\,dr=\frac{2(p-1)}{sp^2}>0.
\]
On $0<t<r<1$, the first argument is $r^s-t^s>0$ and the second is zero. The integrand is strictly positive on this set of positive measure. Since it is nonnegative everywhere else, we conclude \eqref{strictconstant}.

To prove \eqref{cellscale} let us set $r=\ell r'$ and $t=\ell t'$.
Since $\ell>0$, we have $(\ell r')_+=\ell(r')_+$ so
\[
 g(\ell r')=\ell^s g(r')\qquad\text{and}\qquad
 g(\ell r'-\ell)=\ell^s g(r'-1).
\]
The homogeneity of $B_p$ gives a numerator factor $\ell^{sp}$. The denominator contributes
$\ell^{1+sp}$ and $dr\,dt$ contributes $\ell^2$. The resulting factor is
$\ell^{sp+2-(1+sp)}=\ell$. Finally, we prove \eqref{quadratic}. For $p=2$,
\[
 K_{s,2}=\frac12[g-g(\cdot-1)]_{s,2}^2.
\]
With $\widehat v(\xi)=\int_\R e^{-it\xi}v(t)\,dt$, Plancherel gives
\[
 \int_\R|v(x+\tau)-v(x)|^2\,dx
 =\frac1{2\pi}\int_\R|e^{i\tau\xi}-1|^2|\widehat v(\xi)|^2\,d\xi.
\]
Integrate against $|\tau|^{-1-2s}\,d\tau$ and use Tonelli jointly with he identity
\[
 C_{1,s,2}\int_\R\frac{1-\cos(\tau\xi)}{|\tau|^{1+2s}}\,d\tau
 =|\xi|^{2s},
\]
which expresses the Fourier normalization of the fractional Laplacian. Thus
\begin{equation}\label{fractional-plancherel}
 [v]_{s,2}^2
 =\frac1{\pi C_{1,s,2}}\int_\R|\xi|^{2s}|\widehat v(\xi)|^2\,d\xi.
\end{equation}
To apply \eqref{fractional-plancherel} for all $0<s<1$, take
$v_\varepsilon(t)=e^{-\varepsilon t}g(t)-e^{-\varepsilon(t-1)}g(t-1)$, for $\varepsilon>0$. The gamma integral and the translation rule for the Fourier transform give
\[
 \widehat v_\varepsilon(\xi)
 =(1-e^{-i\xi})\Gamma(1+s)(\varepsilon+i\xi)^{-1-s}.
\]
Thus $|\xi|^{2s}|\widehat v_\varepsilon(\xi)|^2$ is dominated by
$\Gamma(1+s)^2|1-e^{-i\xi}|^2/\xi^2$, an integrable function. By dominated convergence, the weighted transforms $|\xi|^s\widehat v_\varepsilon$ are Cauchy in $L^2$. By \eqref{fractional-plancherel} and Fatou's lemma applied to the pointwise limit of $v_\varepsilon$, this gives convergence to $g-g(\cdot-1)$ in the energy seminorm. Applying Plancherel to $\mathbf 1_{(0,1)}$ we conclude
\begin{equation*}
 K_{s,2}
 =\frac{\Gamma(1+s)^2}{2\pi C_{1,s,2}}
  \int_\R\frac{|1-e^{-i\xi}|^2}{\xi^2}\,d\xi
 =\frac{\Gamma(1+s)^2}{C_{1,s,2}} .
\end{equation*}
Substitution in \eqref{cdef} yields \eqref{quadratic}.
\end{proof}

The next lemma records the two regularity estimates used in the proof. The pointwise H\"older bound controls the integrand near the diagonal $x=y$. The local Sobolev bound controls the deformation increments $u_h-u$ in $L^q$.

The interior radius is a fixed small fraction of the distance to the boundary,
\begin{equation}\label{interior-scale}
32r\le d(x_0)\le64r.
\end{equation}
Thus,
\[
 B_{4r}(x_0)\subset B_{8r}(x_0)\subset B_{16r}(x_0)
 \subset B_{24r}(x_0)\subset B_{32r}(x_0)\subset\Omega.
\]
The largest ball is where the equation is used, intermediate balls leave room for interior estimates and cutoffs and the final increment will be estimated on $B_{4r}(x_0)$. More precisely, to estimate the energy remainder for pairs $x,y\in B_r(x_0)$, Lemma \ref{lem:comparison} below requires control of the increment $D_h=u_h-u$ in $B_{4r}(x_0)$. In turn, to control this increment, Lemma \ref{lem:increments} below requires regularity of $u$ in $B_{8r}(x_0)$. Next Lemma \ref{lem:regularity} provides this regularity by using the equation in $B_{32r}(x_0)$.

\begin{samepage}
\begin{lemma}\label{lem:regularity}
Let $u$ be a bounded weak solution of \eqref{eq:dirichlet-problem}. The quotient
$\psi=u/d^s$ extends continuously to $\overline\Omega$ and the zero extension of $u$ belongs to $C^s(\R^N)$. For every
\begin{equation}\label{beta-range}
 s<\beta<\min\{1,sp/(p-1)\},
\end{equation}
one has $|u(x)|\le Cd(x)^s$ and
\begin{equation}\label{weightedholder}
 |u(x)-u(y)|\le Cd(x)^{s-\beta}|x-y|^\beta
 \quad\text{if }x,y\in\Omega,\ |x-y|<d(x)/2.
\end{equation}
Under \eqref{interior-scale}, 
\begin{equation}\label{rescaled}
 U(z)=r^{-s}u(x_0+rz).
\end{equation}
Writing $B_R=B_R(0)$ in these coordinates, one has
$\|U\|_{W^{\alpha,q}(B_8)}\le C_\alpha$ for all $0<\alpha<1$ if
$p\le2$ and for all $0<\alpha<\min\{1,sp/(p-2)\}$ if $p>2$.
The constants are independent of $x_0,r$.
\end{lemma}
\end{samepage}

\begin{proof}
Since $f(u)$ is bounded, the continuity of $u/d^s$ and the bound $|u(x)|\le Cd(x)^s$ follow from \cite[Theorem 1.1]{IMS20,IM24}.
Given $\varphi\in C_c^\infty(B_{32})$, use $\varphi((\cdot-x_0)/r)$, supported in $B_{32r}(x_0)\subset\Omega$, as a test in \eqref{weak}. Taking $x=x_0+rz$ and $y=x_0+r\zeta$ gives
\begin{equation}\label{rescaled-equation}
 (-\Delta)^s_pU=r^sf(r^sU)\qquad\text{weakly in }B_{32}.
\end{equation}
If $x_0+rz\in\Omega$, then
$d(x_0+rz)\le d(x_0)+r|z|\le r(64+|z|)$. Hence
\[
 |U(z)|=r^{-s}|u(x_0+rz)|
 \le C(64+|z|)^s\le C(1+|z|)^s.
\]
If $x_0+rz\notin\Omega$, then $U(z)=0$, so the estimate holds globally.
It gives uniform boundedness on fixed balls. It also controls the exterior values entering the nonlocal equation:
\[
 \int_{|z|>1}\frac{|U(z)|^{p-1}}{|z|^{N+sp}}\,dz
 \le c_1
 \int_1^\infty(1+t)^{s(p-1)}t^{-1-sp}\,dt\le c_2\int_1^\infty t^{-1-s}\,dt=\frac Cs.
\]
The interior estimates \cite[Theorem 1.4]{BLS18}, \cite[Theorem 1.2]{GL24} and \cite[Theorems 1.5 and 1.7]{DKLN25} involve nonlocal tails of the form
\begin{equation}\label{eq:uniform-tail}
 \left(R^{sp}\int_{|z-a|>R}\frac{|U(z)|^{p-1}}{|z-a|^{N+sp}}\,dz\right)^{1/(p-1)}
 \le C .
\end{equation}
The bound is uniform in $x_0,r$ when $a$ ranges over a bounded set and $R$ ranges over a compact subset of $(0,\infty)$. For large $|z|$, $|z-a|$ is comparable to $|z|$ and the preceding integral applies. The same estimate holds with $U$ replaced by its difference from its mean on a fixed interior ball, since those means are uniformly bounded. The right-hand side
$r^sf(r^sU)$ is also uniformly bounded on $B_{32}$, its argument is the bounded original solution and $r$ is bounded in terms of $\Omega$.

Apply \cite[Theorem 1.4]{BLS18} for $p\ge2$ and
\cite[Theorem 1.2]{GL24} for $p<2$ on fixed radius balls inside $B_{32}$.
The local bounds and \eqref{eq:uniform-tail} are uniform on these balls.
A finite cover of $\overline B_{24}$ gives $\|U\|_{C^\beta(B_{24})}\le C$.
Rescaling yields $|u(x)-u(y)|\le Cr^{s-\beta}|x-y|^\beta$ on $B_{24r}(x_0)$. If $|x-y|<d(x)/2$, the segment joining $x,y$ stays at distance comparable to $d(x)$ from the boundary. Then, take $x_0=x$ and $r=d(x)/32$ in \eqref{interior-scale}. Since $|x-y|<16r$, both points lie in
$B_{24r}(x)$ and the preceding estimate gives
\[
 |u(x)-u(y)|\le C r^{s-\beta}|x-y|^\beta \le C d(x)^{s-\beta}|x-y|^\beta.
\]
This proves \eqref{weightedholder}. For $x,y\in\Omega$ with $|x-y|<d(x)/2$, because of \eqref{weightedholder}, we get
\[
 |u(x)-u(y)|\le C|x-y|^s\left(\frac{|x-y|}{d(x)}\right)^{\beta-s}
 \le C|x-y|^s.
\]
If $|x-y|\ge d(x)/2$, then $d(x)\le2|x-y|$ and $d(y)\le d(x)+|x-y|\le3|x-y|$. The estimate  $|u(x)|\le Cd(x)^s$ yields the same bound. If $x\in\Omega$ and $y\notin\Omega$, then $d(x)\le|x-y|$ and the bound follows from $u(y)=0$ and $|u(x)|\le Cd(x)^s$. This proves the global $C^s$ estimate for the zero extension.

The interior $C^\beta(B_{24})$ bound for $U$, with $\beta>s$, provides the initial Sobolev estimate
\[
 [U]_{W^{s,p}(B_{24})}^p
 \le C\int_{B_{24}}\int_{B_{24}}|z-\zeta|^{-N+p(\beta-s)}\,d\zeta\,dz
 \le C\int_0^{48}t^{p(\beta-s)-1}\,dt<\infty.
\]
Thus $\|U\|_{W^{s,p}(B_{24})}\le C$, uniformly in $x_0,r$.
 
Next we translate the regularity of $U$ to regularity of the nonlinearity. To that end, let us apply \cite[Theorem 1.5]{DKLN25} for $p\ge2$ and \cite[Theorem 1.7]{DKLN25} for $1<p<2$, taking their parameter $s_0=s$. Write $p'=p/(p-1)$. If the right-hand side of the equation belongs locally to $W^{t,p'}$, then the solution belongs locally to $W^{\gamma,p}$ for every
\[
 0<\gamma<\min\{1,(sp+t)/(p-1)\}.
\]
Here $t$ is a Sobolev differentiability exponent of the right-hand side and $\gamma$ is a Sobolev differentiability exponent of the solution. The hypotheses on $t$ are
$0<t<\min\{1,sp/(p-2)\}$ for $p>2$ and $0<t<p-1$ for $p\le2$. When $(sp+t)/(p-1)>1$ the results of \cite{DKLN25}  give a weak gradient with positive fractional regularity which implies the weaker assertion $W^{\gamma,p}$ for every $\gamma<1$ used here. 

Next, let $L$ be a Lipschitz constant of $f$ on a bounded interval containing the range of $u$. Then,
\begin{equation}\label{reaction-composition}
 |r^sf(r^sU(z))-r^sf(r^sU(\zeta))|
 \le Lr^{2s}|U(z)-U(\zeta)|.
\end{equation}
The factor $r^{2s}$ is uniformly bounded. On the other hand, on a fixed bounded ball $B$, for $1<a<b$ and $0<t<\gamma<1$, H\"older's inequality gives
\[
 \iint_{B^2}\!\frac{|U(z)-U(\zeta)|^a}{|z-\zeta|^{N+ta}}\,dz\,d\zeta\le
 \left(\iint_{B^2}\!\frac{|U(z)-U(\zeta)|^b}
 {|z-\zeta|^{N+\gamma b}}\,dz\,d\zeta\right)^{a/b}\!\!
 \left(\iint_{B^2}\!|z-\zeta|^{-N+\frac{ab(\gamma-t)}{b-a}}
 \,dz\,d\zeta\right)^{(b-a)/b}.
\]
The second integral is finite since $\gamma>t$. For $a=b$ we use the bounded factor $|z-\zeta|^{(\gamma-t)a}$. The corresponding $L^a$ bound follows from H\"older's inequality on $B$.

For $p\ge2$, take $b=p$ and $a=p'$. Together with \eqref{reaction-composition}, this gives a $W^{t,p'}$ bound for the nonlinearity whenever $0<t<\gamma$.
For $1<p<2$, boundedness of $U$ and \eqref{reaction-composition} instead yield
\[
 |r^sf(r^sU(z))-r^sf(r^sU(\zeta))|^{p'}
 \le C|U(z)-U(\zeta)|^p.
\]
If $t<(p-1)\gamma$, then $tp'<\gamma p$ and
\[
 \frac{|U(z)-U(\zeta)|^p}{|z-\zeta|^{N+tp'}}
 =\frac{|U(z)-U(\zeta)|^p}{|z-\zeta|^{N+\gamma p}}
 |z-\zeta|^{\gamma p-tp'}.
\]
The last factor is bounded on $B^2$, proving the required seminorm bound. The $L^{p'}$ norm of the reaction follows from its uniform boundedness.

Therefore, if $U\in W^{\gamma_k,p}(B)\cap L^\infty(B)$, because of \eqref{reaction-composition}, then $r^sf(r^sU)\in W^{t,p'}(B)$, where $t<\gamma_k$ for $p\ge2$ and  $t<(p-1)\gamma_k$ for $p<2$. This in turn implies (cf. \cite[Theorems 1.5 and 1.7]{DKLN25}) that $U\in W^{\gamma_{k+1},p}(B')$, for $B'\subset B$  and
\[
 \gamma_{k+1}<\min\left\{1,\frac{sp+t}{p-1}\right\}.
\]
An iteration argument shows that the above sequence reaches any $\max\{s,\alpha\}<\gamma_*<1$ for $1<p\leq2$ or $\max\{s,\alpha\}<\gamma_*<\min\{1,sp/(p-2)\}$ when $p>2$ in a finite number of steps. 

If $m$ steps are needed, apply the estimates successively on
\[
 B_{24}\supset B_{24-8/m}\supset\cdots\supset B_{16}.
\]
Each step uses a finite interior covering with radii determined by the gap between consecutive balls. The tail bounds \eqref{eq:uniform-tail} give uniform control in $x_0,r$ on all these balls. Thus $\|U\|_{W^{\gamma_*,p}(B_{16})}\le C$. 

For $p>2$, use the former embedding with exponents $p$ and $2$ and $\gamma_*-\alpha>0$, namely,
\[\frac{|U(z)-U(\zeta)|^2}{|z-\zeta|^{N+2\alpha}}=\left(\frac{|U(z)-U(\zeta)|^p}{|z-\zeta|^{N+p\gamma_*}}\right)^{2/p}|z-\zeta|^{-N(1-2/p)+2(\gamma_*-\alpha)},\]
the first term is bounded in $W^{\gamma_*,p}$ and the second is integrable since $\gamma_*-\alpha>0$. For $p \leq 2$, decrease the order with $q=p$ to obtain the $W^{\alpha,p}$ bound, say, 
\[\frac{|U(z)-U(\zeta)|^p}{|z-\zeta|^{N+p\alpha}}=\frac{|U(z)-U(\zeta)|^p}{|z-\zeta|^{N+p\gamma_*}}|z-\zeta|^{p(\gamma_*-\alpha)},\]
so the seminorm of order $\alpha$ is controlled by that of order $\gamma_*$. Restriction to $B_8$ finishes the proof.
\end{proof}

\begin{remark}
The gap between $B_4$ and $B_8$ is used in Lemma~\ref{lem:increments} below. It keeps the deformed points and all points sampled by mollification inside the ball where the uniform Sobolev estimate holds.
\end{remark}

\section{Comparison with a deformed solution}\label{sec:comparison}

We now explain the deformation argument and then prove the interior estimate that drives the
main proof. There are two steps: first we estimate $D_h=u_h-u$ in $L^q$ using interior regularity; then, by comparison of the equations for $u$ and $u_h$, we estimate the energy remainder, say, the bound \eqref{local-defect} below. On a fixed interior ball, dividing \eqref{local-defect} by $h$ gives $Cr^{N-q\alpha}h^{q\alpha-1}\to0$, since $q\alpha>1$. Thus a fixed interior region leaves no first-order remainder. This argument is not uniform as $r$ approaches the boundary scale $h$; the factor $r^{N-q\alpha}$ records the dependence that must
later be summed over boundary layers.

\medskip
\noindent\textbf{The flow and the associated change of variables.}
The flow of $X$ is the solution of the ordinary differential equation
\begin{equation}\label{flow-ode}
 \frac{d}{dh}\Phi_h(x)=X(\Phi_h(x)),\qquad \Phi_0(x)=x.
\end{equation}
The initial point $x$ labels a trajectory; $h$ is the time parameter.
Since $DX$ is bounded, $X$ is globally Lipschitz and has at most linear growth. The trajectories exist for all small positive and negative times and
\[
 \Phi_{h+k}=\Phi_h\circ\Phi_k,\qquad
 \Phi_h^{-1}=\Phi_{-h},\qquad
 \Phi_h(x)=x+hX(x)+o(h).
\]
The expansion is uniform on fixed bounded sets.
For $X(x)=x$, \eqref{flow-ode} gives $\Phi_h(x)=e^h x$.
Thus $u_h=u\circ\Phi_h$ is a generalization of dilating the argument of $u$.
Initially we require
\begin{equation}\label{outward}
 X\cdot\nu>0\quad\text{on }\partial\Omega.
\end{equation}
Let $\rho$ be signed distance, positive inside $\Omega$. On the boundary,
$\nabla\rho=-\nu$, so along a trajectory
\[
 \left.\frac{d}{dh}\rho(\Phi_h(z))\right|_{h=0}
 =-X(z)\cdot\nu(z)<0,\qquad z\in\partial\Omega.
\]
By compactness this sign holds in a boundary collar. Then, a trajectory starting outside $\Omega$ cannot cross the boundary into $\Omega$: such a crossing would require the signed distance to increase through zero. Hence,
$\Phi_h^{-1}(\Omega)\subset\Omega$ for $h>0$ sufficiently small and
$u_h$ vanishes outside $\Omega$. The flow is uniformly bi-Lipschitz for
these $h$ and the Jacobian determinant of its inverse is uniformly
bounded. As a consequence, 
$\|u_h\|_{L^p(\R^N)}^p\le C\|u\|_{L^p(\R^N)}^p$ and  $[u_h]_{s,p}^p\le C[u]_{s,p}^p$ so $u_h$ and $D_h$ belong to $W^{s,p}_0(\Omega)$.

For pair integrals we use the notation
$\delta v=v(x)-v(y)$ and we also define
\begin{equation}\label{flowkernel}
 \begin{aligned}
 K_h(x,y)&=\frac{C_{N,s,p}}2
 \frac{j_h(x)j_h(y)}{|\Phi_h(x)-\Phi_h(y)|^{N+sp}},\qquad
 K_0(x,y)=\frac{C_{N,s,p}}{2|x-y|^{N+sp}}\\
  j_h(x)&=\det D\Phi_h(x).
 \end{aligned}
\end{equation}
Thus $K_h$ includes both volume Jacobians and the deformed distance.
Differentiating \eqref{flow-ode} in the initial point gives
$\partial_hD\Phi_h=DX(\Phi_h)D\Phi_h$, so 
$\partial_h j_h(x)=j_h(x)\operatorname{tr}(DX(\Phi_h(x))) =j_h(x)\operatorname{div}X(\Phi_h(x))$. Hence,
\[
 j_h(x)=1+h\operatorname{div}X(x)+o(h).
\]
For $x\ne y$, we also have
\[
 \begin{aligned}
 |\Phi_h(x)-\Phi_h(y)|^2
 &=|x-y|^2+2h(X(x)-X(y))\cdot(x-y)+o(h),\\
 \frac{|\Phi_h(x)-\Phi_h(y)|^{-N-sp}}{|x-y|^{-N-sp}}
 &=1-h(N+sp)\frac{(X(x)-X(y))\cdot(x-y)}
 {|x-y|^2}+o(h).
 \end{aligned}
\]
Multiplication by the two Jacobian expansions gives
$K_h=K_0(1+hL_X+o(h))$, with $L_X$ given in \eqref{LX}.
Furthermore, the flow is uniformly bi-Lipschitz:
\begin{equation}\label{biLip}
 e^{-\|DX\|_\infty|h|}|x-y|
 \le|\Phi_h(x)-\Phi_h(y)|
 \le e^{\|DX\|_\infty|h|}|x-y|.
\end{equation}
The analogous bounds for $j_h$, or its differential equation, imply
\begin{equation}\label{kernelbounds}
 |K_h-K_0|\le C|h|K_0,\qquad
 \frac{K_h-K_0}{hK_0}\longrightarrow L_X,\qquad |L_X|\le C.
\end{equation}
Let us stress that the first bound is uniform even though the kernel itself is singular.

We continue with the equation satisfied by $u_h$. If $\varphi$ is supported
inside $\Phi_h^{-1}(\Omega)$, use $\varphi\circ\Phi_{-h}$ as a test in
\eqref{weak} and put the old variables equal to $\Phi_h(x),\Phi_h(y)$.
Then $u(\Phi_h(x))=u_h(x)$ and
$(\varphi\circ\Phi_{-h})(\Phi_h(x))=\varphi(x)$, yielding
\begin{equation}\label{flowweak}
 \iint_{\R^{2N}}J(\delta u_h)\delta\varphi K_h\,dx\,dy
 =\int_{\R^N}j_h f(u_h)\varphi\,dx.
\end{equation}
This forward change of variables will later be distinguished from the inverse change used to differentiate the energy.

Choose $\beta$ as in \eqref{beta-range}. Under \eqref{new-range}, choose
$\alpha$ satisfying
\begin{equation}\label{alpha-range}
 \max\{s,1/q\}<\alpha<\min\{1,s+1/q\},\qquad
 \alpha<\frac{sp}{p-2}\quad\text{if }p>2.
\end{equation}
Lemma~\ref{lem:regularity} gives the uniform $W^{\alpha,q}(B_8)$ bound
for this choice. For $p>2$, the requirement $\alpha>1/2$ is possible
precisely when $sp/(p-2)>1/2$, which explains \eqref{new-range}. For
$p\le2$, the upper Sobolev endpoint is one and \eqref{alpha-range} is
always possible. The layer sums below use the two strict inequalities
\begin{equation}\label{alpha-conditions}
 q\alpha>1,\qquad 1+q(s-\alpha)>0.
\end{equation}
The first gives summability over layers farther from the boundary and
an $o(h)$ bound on each fixed interior ball. The second gives
summability of the increment estimates towards the boundary. The
pointwise exponent $\beta>s$ controls the singularity at $x=y$.

\begin{samepage}
\begin{lemma}\label{lem:increments}
Suppose \eqref{interior-scale} holds and $0<h<cr$ for a fixed $c=c(X)>0$ sufficiently small. Then,
\begin{equation}\label{local-Lq}
 \int_{B_{4r}(x_0)}|D_h|^q\,dx
 \le Ch^{q\alpha}r^{N+q(s-\alpha)}
\end{equation}
and
\begin{equation}\label{ball-Lq}
 \int_{B_L(x_0)}|D_h|^q\,dx
 \le Ch^{q\alpha}L^{N+q(s-\alpha)}
 \qquad (L\ge4r).
\end{equation}
Moreover, for $t\ge C_0h$ with $C_0>32/c$ fixed,
\begin{equation}\label{layer-Lq}
 \int_{\{x\in\Omega:t\le d(x)<2t\}}|D_h|^q\,dx
 \le Ch^{q\alpha}t^{1+q(s-\alpha)}.
\end{equation}
The thin-layer $d<C_0h$ satisfies the estimate
\[
 \int_{\{x\in\Omega:d(x)<C_0h\}}|D_h(x)|^q\,dx\le Ch^{sq+1},
\]
and
\begin{equation}\label{global-Lq}
 \|D_h\|_{L^q(\Omega)}^q\le Ch^{q\alpha}.
\end{equation}
\end{lemma}
\end{samepage}
\begin{proof}
In the coordinates $x=x_0+rz$,
\[
 D_h(x_0+rz)=r^s
 \left[U\left(\frac{\Phi_h(x_0+rz)-x_0}{r}\right)-U(z)\right].
\]
Since $\Phi_h(x)-x=\int_0^hX(\Phi_\tau(x))\,d\tau$, boundedness of $X$ on a fixed neighborhood gives $|\Phi_h(x)-x|\le Ch$. The rescaled displacement is therefore at most $Ch/r$ and
\[
 D_z\left(\frac{\Phi_h(x_0+rz)-x_0}{r}\right)=D\Phi_h(x_0+rz)=I+O(h).
\]
Choose $c$ so that $Ch/r<1$ whenever $h<cr$. For $z\in B_4$, each point
of the segment joining $z$ to $(\Phi_h(x_0+rz)-x_0)/r$ is within distance
$Ch/r$ of $z$ and hence lies in $B_5\subset B_6$.

Let $U_\varepsilon=\eta_\varepsilon*U$, where $\varepsilon=h/r<1$ and $\eta_\varepsilon$ is a nonnegative smooth mollifier supported in $B_\varepsilon$, of integral one. Jensen's inequality give
\[
\begin{aligned}
 \|U-U_\varepsilon\|_{L^q(B_6)}^q
 &\le C\varepsilon^{-N}\int_{|w|<\varepsilon}\int_{B_6}|U(z-w)-U(z)|^q\,dz\,dw,\\
 \|\nabla U_\varepsilon\|_{L^q(B_6)}^q
 &\le C\varepsilon^{-N-q}\int_{|w|<\varepsilon}\int_{B_6}|U(z-w)-U(z)|^q\,dz\,dw.
\end{aligned}
\]
Inserting the denominator $|w|^{N+\alpha q}$ and using $|w|<\varepsilon$ yields
\begin{equation}\label{mollifier-bounds}
 \|U-U_\varepsilon\|_{L^q(B_6)}\le C\varepsilon^\alpha\|U\|_{W^{\alpha,q}(B_8)}\qquad\text{and}\qquad
 \|\nabla U_\varepsilon\|_{L^q(B_6)}\le C\varepsilon^{\alpha-1}\|U\|_{W^{\alpha,q}(B_8)}.
\end{equation}
Next, let us write
\[
\begin{aligned}
U\left(\frac{\Phi_h(x_0+rz)-x_0}{r}\right)-U(z)&= (U-U_\varepsilon)\left(\frac{\Phi_h(x_0+rz)-x_0}{r}\right)
 +(U_\varepsilon-U)(z)\\
 &\quad+U_\varepsilon\left(\frac{\Phi_h(x_0+rz)-x_0}{r}\right)-U_\varepsilon(z).
\end{aligned}
\]
By \eqref{mollifier-bounds} and Lemma~\ref{lem:regularity}, $\|U_\varepsilon-U\|_{L^q(B_4)}\le C\varepsilon^\alpha$. Changing variables in the error at the deformed point gives
\[
 \int_{B_4}\left|(U-U_\varepsilon)\left(\frac{\Phi_h(x_0+rz)-x_0}{r}\right)\right|^qdz
 \le C\int_{B_6}|U-U_\varepsilon|^q\le C\varepsilon^{\alpha q}.
\]
For the smooth term, integrate the gradient along the joining segment:
{\small \[
\begin{aligned}
U_\varepsilon\left(\frac{\Phi_h(x_0+rz)-x_0}{r}\right)-U_\varepsilon(z)=\int_0^1\nabla U_\varepsilon\left(z+\theta\frac{\Phi_h(x_0+rz)-(x_0+rz)}r\right)
 \cdot\frac{\Phi_h(x_0+rz)-(x_0+rz)}r\,d\theta.
\end{aligned}
\]}
For each $\theta\in[0,1]$, the argument of the gradient defines a map
from $B_4$ into $B_6$ with derivative $I+\theta(D\Phi_h-I)=I+O(h)$.
The derivative bound and convexity of $B_4$ give injectivity and the
inverse Jacobian is bounded uniformly in $\theta$. Applying Jensen's
inequality in $\theta$ and then changing variables in $z$ gives
\[
\int_{B_4}\left|U_\varepsilon\left(\frac{\Phi_h(x_0+rz)-x_0}{r}\right)-U_\varepsilon(z)\right|^q\,dz
 \le C(h/r)^q\|\nabla U_\varepsilon\|_{L^q(B_6)}^q
 \le C(h/r)^{q\alpha},
\]
where the last step uses \eqref{mollifier-bounds} and the uniform bound
from Lemma~\ref{lem:regularity}. Gathering together the former bounds, we conclude
\[
 \int_{B_{4r}(x_0)}|D_h|^qdx=r^{N+sq}\int_{B_4}\left|U\left(\frac{\Phi_h(x_0+rz)-x_0}{r}\right)-U(z)\right|^qdz\le Ch^{q\alpha}r^{N+q(s-\alpha)}.
\]
This proves \eqref{local-Lq}. Next we prove \eqref{layer-Lq}. For small $t$, tubular coordinates give
\[
 |\{x\in\Omega:t\le d(x)<2t\}|
 =\int_{\partial\Omega}\int_t^{2t}(1+O(\tau))\,d\tau\,d\sigma(z)\le Ct.
\]
The same bound holds for the enlarged layer $t/2<d<3t$. Choose a maximal collection of centers in $t\le d<2t$ with mutual
distances at least $t/32$, i.e., the balls of radius $t/32$ cover the layer.
The balls of radius $t/64$ are disjoint and the $1$-Lipschitz
property of the distance puts them in $t/2<d<3t$ inside $\Omega$.
Each such ball has volume of order $t^N$, so their number is at most
$Ct/t^N=Ct^{1-N}$. At every center the radius $t/32$ satisfies
\eqref{interior-scale} and $t\ge C_0h$ ensures $h<ct/32$ because
$C_0>32/c$. Applying \eqref{local-Lq} on each ball and summing we conclude \eqref{layer-Lq}, say,
\[
 \int_{\{x\in\Omega:t\le d(x)<2t\}}|D_h|^q\,dx
 \le Ct^{1-N}\,h^{q\alpha}t^{N+q(s-\alpha)}
 =Ch^{q\alpha}t^{1+q(s-\alpha)}.
\]
By Lemma~\ref{lem:regularity}, the zero extension of $u$ is in
$C^s(\R^N)$. Thus for $x\in\Omega$,
\[
 |D_h(x)|=|u(\Phi_h(x))-u(x)|\le C|\Phi_h(x)-x|^s\le Ch^s.
\]
For small $h$, tubular coordinates give the volume bound $Ch$ on $d<C_0h$. For larger $h$, the bound follows from the boundedness of $\Omega$. Hence, this part contributes at most $Ch^{sq+1}$. For the rest of $\Omega$, apply
\eqref{layer-Lq} at $t=2^kC_0h$ with $k\in\mathcal{K}=\{k\ge0: 2^kC_0h\le\operatorname{diam}(\Omega)\}$. These layers cover all depths
$d\ge C_0h$ because $d\le\operatorname{diam}(\Omega)<\infty$. Using \eqref{alpha-conditions}, we obtain
\[
 \|D_h\|_{L^q(\Omega)}^q
 \le Ch^{sq+1}
 +Ch^{q\alpha}\sum_{k\in\mathcal{K}}
 \le Ch^{q\alpha}.
\]
This proves \eqref{global-Lq}. To sum the same local estimate inside $B_L(x_0)$, say  \eqref{ball-Lq}, use the volume bound
\begin{equation}\label{collar-volume}
 |\Omega\cap B_L(x_0)\cap\{d<t\}|\le CL^{N-1}t.
\end{equation}
When $t\le L$ is within the tubular radius, the nearest boundary
points lie in $B_{2L}(x_0)$ and boundary graphs give
$\sigma(\partial\Omega\cap B_{2L}(x_0))\le CL^{N-1}$. Integrating in
the normal direction gives \eqref{collar-volume}. For $t>L$ the
volume of $B_L$ suffices. If $t\le L$ and $t$ is bounded below by the
tubular radius, the bound follows from the volume of $\Omega$. 

Fix $t\ge C_0h$ and consider the portion of the layer
$\{t\le d<2t\}$ inside $B_L(x_0)\cap\Omega$. By
\eqref{interior-scale} and $L\ge 4r$, every point of this set satisfies
\[
 t\le d(x)\le d(x_0)+|x-x_0|
 \le 64r+L\le 17L.
\]
Choose a maximal collection of centers $x_i$ in this set with mutual
distances at least $t/32$ so the balls $B_{t/32}(x_i)$ cover the set,
and the balls $B_{t/64}(x_i)$ are pairwise disjoint. Observe that $t/64<d(x_i)$, so $B_{t/64}(x_i)\subset\Omega$.
Moreover, for every $y\in B_{t/64}(x_i)$, we have $|y-x_0|<L+t/64\le L+17L/64<2L$,
 and $d(y)\le d(x_i)+|y-x_i|<2t+t/64<3t$. Thus \eqref{collar-volume}, applied with $2L$ and $3t$, gives
\[
 \sum_i |B_{t/64}(x_i)|
 \le |\Omega\cap B_{2L}(x_0)\cap\{d<3t\}|
 \le CL^{N-1}t.
\]
Since each ball has volume $\omega_N(t/64)^N$, there are at most
$CL^{N-1}t^{1-N}$ centers.
At each center, the radius $t/32$ satisfies \eqref{interior-scale}
because $t\le d(x_i)<2t$. Also, $t\ge C_0h$ and $C_0>32/c$ imply
$h<c(t/32)$. We may therefore apply \eqref{local-Lq} with center
$x_i$ and radius $t/32$. Restricting that estimate to the covering
ball $B_{t/32}(x_i)$ and summing over the centers, we obtain
\[
 \begin{aligned}
 \int\limits_{B_L(x_0)\cap\Omega\cap\{t\le d<2t\}}\mkern-7mu |D_h|^q\,dx
 \le \sum_i\int\limits_{B_{t/32}(x_i)}\mkern-5mu |D_h|^q\,dx&\le CL^{N-1}t^{1-N}\,h^{q\alpha}t^{N+q(s-\alpha)}=Ch^{q\alpha}L^{N-1}t^{1+q(s-\alpha)}.
 \end{aligned}
\]
On the remaining layer $\{d<C_0h\}$, the pointwise estimate
$|D_h|\le Ch^s$ and \eqref{collar-volume} yield
\[
 \int\limits_{B_L(x_0)\cap\Omega\cap\{d<C_0h\}}\mkern-7mu |D_h|^q\,dx
 \le Ch^{sq}|B_L(x_0)\cap\Omega\cap\{d<C_0h\}|
 \le CL^{N-1}h^{sq+1}.
\]
Since $D_h=0$ outside $\Omega$, it remains to sum the preceding
layer estimates at the scales $t=2^jC_0h$ with $j\in\mathcal{J}=\{j\geq0: 2^jC_0h\le 17L\}$. Hence,
\[
 \int_{B_L(x_0)} |D_h|^q\,dx
 \le CL^{N-1}h^{sq+1}
 +Ch^{q\alpha}L^{N-1}
 \sum_{j\in\mathcal{J}}
 (2^jC_0h)^{1+q(s-\alpha)}\le Ch^{q\alpha}L^{N+q(s-\alpha)},
\]
since $1+q(s-\alpha)>0$ and $h\le L$ so
$h^{sq+1}=h^{q\alpha}h^{1+q(s-\alpha)}
\le h^{q\alpha}L^{1+q(s-\alpha)}$. 
\end{proof}

The next lemma is a local energy estimate obtained from the equation.

\begin{samepage}

\begin{lemma}\label{lem:comparison}
Under \eqref{outward}, \eqref{beta-range} and \eqref{alpha-range},
\begin{equation}\label{local-defect}
 \iint_{B_r(x_0)^2}B_p(\delta u,\delta u_h)K_0\,dx\,dy
 \le Cr^N(h/r)^{q\alpha}
\end{equation}
whenever \eqref{interior-scale} holds and $0<h<cr$.
\end{lemma}
\end{samepage}
\begin{proof}
All balls in the proof are centered at $x_0$. Choose a smooth cutoff
$0\le\chi\le1$ with $\chi=1$ on $B_r$, $\chi=0$ off $B_{2r}$ and
$|\nabla\chi|\le C/r$. For $h/r$ sufficiently small, the flow estimate gives
\[
 |\Phi_h(x)-x_0|\le |x-x_0|+Ch<3r
 \quad (x\in\overline B_{2r}).
\]
By \eqref{interior-scale}, $\overline B_{3r}\subset\Omega$, so
$\overline B_{2r}\subset\Omega\cap\Phi_h^{-1}(\Omega)$.
Multiplication by the bounded Lipschitz function $\chi^q$ preserves the
$W^{s,p}(\R^N)$ regularity and
$\varphi=\chi^qD_h$ is admissible in both \eqref{weak} and
\eqref{flowweak}. Subtracting the equations and using
\[
 J(\delta u_h)K_h-J(\delta u)K_0
 =[J(\delta u_h)-J(\delta u)]K_0+J(\delta u_h)(K_h-K_0),
\]
gives
\begin{equation}\label{comparison-exact}
\begin{split}
 \iint_{\R^{2N}}
 [J(\delta u_h)-J(\delta u)]\delta(\chi^qD_h)K_0\,dx\,dy=&\int_{\R^N}[j_hf(u_h)-f(u)]\chi^qD_h\,dx\\
 & +\iint_{\R^{2N}}(K_0-K_h)J(\delta u_h)
 \delta(\chi^qD_h)\,dx\,dy.
\end{split}
\end{equation}
The pair integrals are finite by the energy bounds and H\"older's inequality. Set
\[
 S(x,y)=[J(\delta u_h)-J(\delta u)]\,\delta D_h.
\]
Note that, by \eqref{Sidentity}, $S=B_p(\delta u,\delta u_h)+B_p(\delta u_h,\delta u) \ge B_p(\delta u,\delta u_h)$. We use the product identity
\begin{equation}\label{product-id}
 \delta(\chi^qD_h)
 =\frac{\chi(x)^q+\chi(y)^q}{2}\delta D_h
 +\frac{D_h(x)+D_h(y)}2\,\delta(\chi^q).
\end{equation}
The integrand on the left of \eqref{comparison-exact} is symmetric in
$x,y$. It vanishes on $(B_{4r}^c)^2$ and on $B_{4r}\times B_{4r}^c$ the
test difference is $\chi(x)^qD_h(x)$, supported in $x\in B_{2r}$.
Consequently,
\begin{align}
 \iint_{\R^{2N}}[J(\delta u_h)-J(\delta u)]
 \delta(\chi^qD_h)K_0\,dx\,dy&=\iint_{B_{4r}^2}\frac{\chi(x)^q+\chi(y)^q}{2}S K_0\,dx\,dy\nonumber\\
 &\quad+\iint_{B_{4r}^2}[J(\delta u_h)-J(\delta u)]
 \frac{D_h(x)+D_h(y)}2\delta(\chi^q)K_0\,dx\,dy\nonumber\\
 &\quad+2\int_{B_{2r}}\int_{B_{4r}^c}[J(\delta u_h)-J(\delta u)]
 \chi(x)^qD_h(x)K_0\,dy\,dx.\label{comparison-split}
\end{align}
We use the equality in \eqref{comparison-split} above to obtain a bound for the term
\begin{equation}\label{principal-term}
 \iint_{B_{4r}^2}\frac{\chi(x)^q+\chi(y)^q}{2}S(x,y)K_0\,dx\,dy.
\end{equation}
Since $\chi=1$ on $B_r$ and  $S \ge B_p(\delta u,\delta u_h)$ this will produce the bound \eqref{local-defect}.

For $p\ge2$ we have $q=2$. Combining \eqref{Munified} and
\eqref{weighted-coercivity} yields
\[
 |J(\delta u_h)-J(\delta u)|
 \le CS^{1/2}(|\delta u|+|\delta u_h|)^{(p-2)/2}.
\]
Since $\delta(\chi^2)=(\chi(x)+\chi(y))\delta\chi$,
the absolute value of the integrand of the second term in
\eqref{comparison-split}, before multiplication by $K_0$, is at most
\[
 C\left(\frac{\chi(x)^2+\chi(y)^2}{2}S\right)^{1/2}
 (|\delta u|+|\delta u_h|)^{(p-2)/2}
 (|D_h(x)|+|D_h(y)|)|\delta\chi|.
\]
Then, for any $\varepsilon>0$, Young's inequality bounds the above quantity by
\begin{equation}\label{eq:pp}
 \varepsilon\frac{\chi(x)^2+\chi(y)^2}{2}S
 +C_\varepsilon (|\delta u|+|\delta u_h|)^{p-2}
 (|D_h(x)|^2+|D_h(y)|^2)|\delta\chi|^2.
\end{equation}

For $x,y\in B_{4r}$, both points and their images under $\Phi_h$ belong
to $B_{5r}$. By \eqref{interior-scale}, their boundary distances are
comparable to $r$ and their pairwise distances are smaller than half
the corresponding boundary distance. Thus \eqref{weightedholder} and
the Lipschitz bound on $\Phi_h$ give
\begin{equation}\label{Hö}
 |\delta u|+|\delta u_h|\le Cr^{s-\beta}|x-y|^\beta.
\end{equation}

Integration against $K_0$ in \eqref{eq:pp} gives $\varepsilon$ times
\eqref{principal-term} and, using $|\delta\chi|\le C|x-y|/r$ and \eqref{local-Lq}, 
\begin{equation*}
\begin{split}
\iint_{B_{4r}^2} (|\delta u|+|\delta u_h|)^{p-2} &(|D_h(x)|^2+|D_h(y)|^2)|\delta\chi|^2K_0\,dx\,dy\\
 &\le c_1 r^{(s-\beta)(p-2)-2}\int_{B_{4r}}|D_h(x)|^2\int_{B_{4r}} |x-y|^{\beta(p-2)+2-N-sp}\,dy dx\\
&\le c_3 r^{-2s}\int_{B_{4r}} |D_h(x)|^2\,dx\\
&
\le
c_4 r^N\left(\frac{h}{r}\right)^{2\alpha},
\end{split}
\end{equation*}
For $1<p<2$ we have $q=p$ and \eqref{Msmall} gives
$|J(\delta u_h)-J(\delta u)|\le CS^{(p-1)/p}$. The mean value theorem jointly with $\max\{\chi(x),\chi(y)\}^p\le \chi(x)^p+\chi(y)^p$ give
\[
 |\delta(\chi^p)|\le p(\chi(x)^p+\chi(y)^p)^{(p-1)/p}|\delta\chi|.
\]
Multiplying these bounds by $|D_h(x)+D_h(y)|/2$ and using Young inequality produces the bound
\[
 \varepsilon\frac{\chi(x)^p+\chi(y)^p}{2}S
 +C_\varepsilon(|D_h(x)|^p+|D_h(y)|^p)|\delta\chi|^p.
\]
The integral of the first summand against $K_0$ is again
$\varepsilon$ times \eqref{principal-term} and, by \eqref{local-Lq},
\[\begin{split}
\iint_{B_{4r}^2}
 (|D_h(x)|^p+|D_h(y)|^p)|\delta\chi|^p K_0\,dx\,dy&\le  r^{-p}\int_{B_{4r}}|D_h(x)|^p
 \left(\int_{B_{4r}}|x-y|^{p-N-sp}\,dy\right) dx\\
 &\le c_1 r^N\left(\frac{h}{r}\right)^{q\alpha}.
\end{split}
\]
For $x\in B_{2r}$ and $y\in B_{4r}^c$ one has $|x-y|\ge2r$. The global
$C^s$ bound in Lemma~\ref{lem:regularity} and the uniform Lipschitz
bound on $\Phi_h$ imply
\[
 |\delta u|+|\delta u_h|\le C|x-y|^s.
\]
Thus, by \eqref{Munified} and $\delta D_h=\delta u_h-\delta u$,
we get following pointwise estimate for the last term of
\eqref{comparison-split},
\[
 \begin{aligned}
 |J(\delta u_h)-J(\delta u)|\,\chi(x)^q|D_h(x)|K_0&\le
 C|x-y|^{s(p-q)-N-sp}|\delta D_h|^{q-1}|D_h(x)|\\
 &\le
 C\frac{|D_h(x)|^q+|D_h(x)||D_h(y)|^{q-1}}{|x-y|^{N+sq}}.
 \end{aligned}
\]
The integral of the first summand is controlled by \eqref{local-Lq},
since
\[
 \int_{B_{2r}}\int_{B_{4r}^c}
 \frac{|D_h(x)|^q}{|x-y|^{N+sq}}\,dy\,dx
 \le \frac{c}{sq(2r)^{sq}}
 \int_{B_{2r}}|D_h(x)|^q\,dx
 \le Cr^N\left(\frac{h}{r}\right)^{q\alpha}.
\]
Next, for each fixed $x\in B_{2r}$, take $B_{4r}^c=\cup_{k\geq2}\mathcal{A}_k$ with $\mathcal{A}_k=\{2^kr\le|y-x_0|<2^{k+1}r\}$. Then
\[
 \begin{aligned}
 \int_{B_{4r}^c}\frac{|D_h(y)|^{q-1}}{|x-y|^{N+sq}}\,dy
 &=\sum_{k=2}^{\infty}\int_{\mathcal{A}_k}
 \frac{|D_h(y)|^{q-1}}{|x-y|^{N+sq}}\,dy\le C\sum_{k=2}^{\infty}(2^kr)^{-N-sq}
 \int_{B_{2^{k+1}r}} |D_h(y)|^{q-1}\,dy .
 \end{aligned}
\]
Here $|x-y|\ge |y-x_0|-2r\ge 2^{k-1}r$ on $\mathcal{A}_k$.
H\"older's inequality and \eqref{ball-Lq} give
\[
 \int_{B_L}|D_h|^{q-1}\le |B_L|^{1/q}\left(\int_{B_L}|D_h|^q\right)^{(q-1)/q}\le Ch^{\alpha(q-1)}L^{N+(s-\alpha)(q-1)},
 \qquad \text{for $L\ge4r$}.
\]
With $L=2^{k+1}r$ this controls the $k$th summand by
$Ch^{\alpha(q-1)}(2^kr)^{-s-\alpha(q-1)}$. Summing the geometric series
and using H\"older's inequality with \eqref{local-Lq} on $B_{2r}$ gives
\[
 \begin{aligned}
 \int_{B_{2r}}\int_{B_{4r}^c}
 \frac{|D_h(x)||D_h(y)|^{q-1}}{|x-y|^{N+sq}}\,dy\,dx
 &\le Ch^{\alpha(q-1)}r^{-s-\alpha(q-1)}
 |B_{2r}|^{1-1/q}\left(\int_{B_{2r}}|D_h|^q\right)^{1/q}\\
 &\le Ch^{\alpha(q-1)}r^{-s-\alpha(q-1)}
 h^\alpha r^{N+s-\alpha}\\
 &=Cr^N\left(\frac{h}{r}\right)^{q\alpha}.
 \end{aligned}
\]

Next, we address the change of kernel. Using \eqref{product-id}, we get
\begin{align*}
 \iint_{\R^{2N}}(K_0-K_h)J(\delta u_h) \delta(\chi^qD_h)\,dx\,dy
 &=\iint_{B_{4r}^2}(K_0-K_h)J(\delta u_h)
 \frac{\chi(x)^q+\chi(y)^q}{2}\delta D_h\,dx\,dy\\
 &\quad+\iint_{B_{4r}^2}(K_0-K_h)J(\delta u_h)
 \frac{D_h(x)+D_h(y)}2\,\delta(\chi^q)\,dx\,dy\\
 &\quad+2\int_{B_{2r}}\int_{B_{4r}^c}
 (K_0-K_h)J(\delta u_h)\chi(x)^qD_h(x)\,dy\,dx .
\end{align*}
For the first summand, \eqref{weighted-coercivity} yields
\[
 |J(\delta u_h)\delta D_h|
 \le CS^{1/2}(|\delta u|+|\delta u_h|)^{p/2}.
\]
Using
$|K_h-K_0|\le ChK_0$ from \eqref{kernelbounds} and Young's inequality gives
\[(K_0-K_h)J(\delta u_h)
 \frac{\chi(x)^q+\chi(y)^q}{2}\delta D_h\leq  \varepsilon\frac{\chi(x)^q+\chi(y)^q}{2}S K_0
 +C_\varepsilon h^2 (|\delta u|+|\delta u_h|)^pK_0 .\]

The first term in the right hand side integrates to $\varepsilon$ times
\eqref{principal-term} while, by using \eqref{Hö} and $\beta>s$,
\[
 C_\varepsilon h^2 \iint_{B_{4r}^2} (|\delta u|+|\delta u_h|)^pK_0\,dx\,dy\le Ch^2r^{p(s-\beta)}|B_{4r}|\int_0^{8r}\tau^{p\beta-sp-1}\,d\tau\le Ch^2r^N,
\]
For the second local summand, use
$|\delta(\chi^q)|\le C|\delta\chi|\le C|x-y|/r$.
The kernel estimate and the same local H\"older bound imply
\[
 \left|(K_0-K_h)J(\delta u_h)
 \frac{D_h(x)+D_h(y)}2\,\delta(\chi^q)\right|\le Ch r^{(s-\beta)(p-1)-1}
 \frac{|D_h(x)|+|D_h(y)|}{|x-y|^{N+sp-\beta(p-1)-1}}.
\]
By symmetry its integral over $B_{4r}^2$ is at most
\[
 Ch r^{(s-\beta)(p-1)-1}
 \int_{B_{4r}} |D_h(x)|
 \int_{B_{4r}}|x-y|^{\beta(p-1)+1-N-sp}\,dy dx
\le Ch r^{-s}\int_{B_{4r}}|D_h|
 \le Cr^Nh\left(\frac{h}{r}\right)^\alpha.
\]
The last inequality follows from \eqref{local-Lq} and Hölder inequality. 

For the cross term, $|J(\delta u_h)|\le C|x-y|^{s(p-1)}$ by global
$C^s$ regularity and the Lipschitz bound for the flow. Its absolute value is at most
\[
 Ch\int_{B_{2r}}|D_h(x)|\int_{B_{4r}^c}|x-y|^{-N-s}\,dy\,dx\le Ch r^{-s}\int_{B_{2r}}|D_h|\le Cr^Nh\left(\frac{h}{r}\right)^\alpha .
\]
Thus, 
\begin{equation}\label{kernel-error} 
\begin{split}
\iint_{\R^{2N}}(K_0-K_h)J(\delta u_h) \delta(\chi^qD_h)\,dx\,dy\leq &\varepsilon  \iint_{B_{4r}^2}\frac{\chi(x)^q+\chi(y)^q}{2}S(x,y)K_0\,dx\,dy\\
&+C_\varepsilon r^N\left(h^2+h\left(\frac{h}{r}\right)^\alpha\right).
\end{split}
\end{equation}

Finally, we address the term involving the nonlinearity $f$. First, note that, on the support of $\chi$, 
\[
 \int_{\R^N}[j_hf(u_h)-f(u)]\chi^qD_h\,dx=\int_{B_{2r}}j_h[f(u_h)-f(u)]\chi^qD_h\,dx
 +\int_{B_{2r}}(j_h-1)f(u)\chi^qD_h\,dx .
\]
Since $u$ and $u_h$ are bounded, $f$ is Lipschitz on an interval containing their values and
$j_h=1+O(h)$, the two integrals are bounded in absolute value by
$C\int_{B_{2r}}|D_h|^2$ and $Ch\int_{B_{2r}}|D_h|$, respectively.
The growth bound in \eqref{weightedholder} gives $|D_h|\le Cr^s$ on
$B_{2r}$, as both $x$ and $\Phi_h(x)$ lie in $B_{3r}$.
Since $q\le2$, \eqref{local-Lq} then yields
\[
 \begin{aligned}
 \int_{B_{2r}} |D_h|^2
 &\le Cr^{s(2-q)}\int_{B_{2r}}|D_h|^q
 \le Cr^{N+2s}\left(\frac{h}{r}\right)^{q\alpha},\\
 h\int_{B_{2r}}|D_h|
 &\le h|B_{2r}|^{1-1/q}\left(\int_{B_{2r}}|D_h|^q\right)^{1/q}
 \le Cr^{N+s}h\left(\frac{h}{r}\right)^\alpha.
 \end{aligned}
\]
As $r$ is bounded, their sum is at most
$Cr^N((h/r)^{q\alpha}+h(h/r)^\alpha)$. Combining \eqref{comparison-exact} and \eqref{comparison-split} with
the preceding estimates gives
\begin{align*}
 \iint_{B_{4r}^2}\frac{\chi(x)^q+\chi(y)^q}{2}S(x,y)K_0\,dx\,dy\le&2\varepsilon\iint_{B_{4r}^2}\frac{\chi(x)^q+\chi(y)^q}{2}S(x,y)K_0\,dx\,dy\\
 &+C_\varepsilon r^N\left[\left(\frac hr\right)^{q\alpha}+h^2+h\left(\frac hr\right)^\alpha\right].
\end{align*}
Choose $\varepsilon>0$ so that $2\varepsilon\le 1/2$ and move
the first term on the right to the left. Since $S\ge0$ and
$\chi=1$ on $B_r$, we obtain
\begin{align*}
 \iint_{B_r^2} S K_0\,dx\,dy
 &\le \iint_{B_{4r}^2}\frac{\chi(x)^q+\chi(y)^q}{2}
       S K_0\,dx\,dy\le Cr^N\left[\left(\frac hr\right)^{q\alpha}
       +h^2+h\left(\frac hr\right)^\alpha\right].
\end{align*}
Since $0<h/r<1$, $r$ is bounded and
$q\alpha<2$ and $q\alpha<1+\alpha$, we have
\[
 h^2=r^2\left(\frac hr\right)^2
 \le C\left(\frac hr\right)^{q\alpha},
 \qquad
 h\left(\frac hr\right)^\alpha
 =r\left(\frac hr\right)^{1+\alpha}
 \le C\left(\frac hr\right)^{q\alpha}.
\]
Consequently, by \eqref{Sidentity} and the nonnegativity of $B_p$, we conclude 
\[
 \iint_{B_r^2} B_p(\delta u,\delta u_h)K_0\,dx\,dy
 \le \iint_{B_r^2} S K_0\,dx\,dy
 \le Cr^N\left(\frac hr\right)^{q\alpha}.
\]
\end{proof}
The local comparison relies on the bounds for $D_h$ jointly with the
interior H\"older estimate with $\beta>s$, which makes the terms near
the diagonal integrable. The two conditions on $\alpha$ enter in
different sums. The inequality $1+q(s-\alpha)>0$ allows us to sum the
increments over the layers in a large ball, down to the boundary.
The inequality $q\alpha>1$ will make the contribution to the energy
remainder, divided by $h$, small for interior pairs at depths
$t\ge Rh$ as $R\to\infty$ (see \eqref{near_ine} below). For $p>2$, the bound
$\alpha<sp/(p-2)$ comes from the Sobolev estimate in
Lemma~\ref{lem:regularity} used to control $D_h$.

\section{Boundary concentration and the Pohozaev identity}\label{sec:boundary}

To compute the first-order contribution to $R_h$ in \eqref{Rdef}, we first
keep $x$ within distance $Rh$ of the boundary and restrict the separation
to $\varepsilon h<|x-y|<Rh$. Rescaling gives the limit on this set and
hence a lower bound. For the upper bound, we estimate the omitted pairs:
the H\"older bounds control small separations in the boundary strip,
Lemma~\ref{lem:comparison} controls nearby pairs in the interior and
Lemma~\ref{lem:increments} controls separated pairs.

\begin{samepage}
\begin{proposition}\label{prop:boundary}
Let $0<s<1$ and assume that $X$ satisfies \eqref{outward}. If
$u\in W^{s,p}_0(\Omega)$ and $u/d^s$ extends continuously to
$\overline\Omega$, then, 
\begin{equation}\label{lower}
 \liminf_{h\downarrow0}\frac{R_h}{h}
 \ge c_{s,p}\int_{\partial\Omega}|\psi|^pX\cdot\nu\,d\sigma.
\end{equation}
Under the hypotheses of Theorem~\ref{thm:general}, the limit exists and
\begin{equation}\label{flowlimit}
 \lim_{h\downarrow0}\frac{R_h}{h}
 =c_{s,p}\int_{\partial\Omega}|\psi|^pX\cdot\nu\,d\sigma.
\end{equation}
\end{proposition}
\end{samepage}
\begin{proof}
Let $\rho$ be the signed distance to $\partial\Omega$, positive inside
and negative outside. It is $C^2$ on a fixed tubular neighborhood, where
every point has a unique representation
$x=z-\rho(x)\nu(z)$, $z\in\partial\Omega$. Its volume element is
$(1+O(|\rho|))\,d\rho\,d\sigma(z)$. First fix $R>1$ and $0<\varepsilon<R$. Retain only pairs in
\begin{equation}\label{core}
 \mathcal C_{h,R,\varepsilon}
 =\{(x,y):|\rho(x)|<Rh,\ \varepsilon h<|x-y|<Rh\}.
\end{equation}

Partitioning the pairs $(x,y)\in\R^{2N}$ off the diagonal gives
\begin{equation}\label{pair-decomposition}
\begin{aligned}
R_h=&\iint_{\mathcal C_{h,R,\varepsilon}}B_p(\delta u,\delta u_h)K_0\,dx\,dy+\iint_{\substack{|\rho(x)|<Rh\\0<|x-y|\le\varepsilon h}}B_p(\delta u,\delta u_h)K_0\,dx\,dy\\
&+\iint_{\{|x-y|\ge Rh\}}B_p(\delta u,\delta u_h)K_0\,dx\,dy+\iint_{\substack{|\rho(x)|\ge Rh\\0<|x-y|<Rh}}B_p(\delta u,\delta u_h)K_0\,dx\,dy.
\end{aligned}
\end{equation}
We rescale the first integral below and evaluate the resulting profile in
Step~2. Step~3 estimates the second integral as $\varepsilon\downarrow0$.
Steps~4--6 estimate the last two integrals together as $R\to\infty$,
according to how $|x-y|$ compares with the distances to the boundary.

The upper cutoffs keep all rescaled coordinates bounded. Let
$z\in\partial\Omega$ be the nearest boundary point to $x$ and write
\[
 x=z-hr\nu(z),\qquad y=x+hw.
\]
Plainly $\rho(x)=hr$ and $|x-y|=h|w|$, so the conditions defining
$\mathcal C_{h,R,\varepsilon}$ become $|r|<R$ and
$\varepsilon<|w|<R$. The lower cutoff keeps the rescaled separation
away from the singularity at $w=0$.

Set $t=r-w\cdot\nu(z)$. Since $\rho(z)=0$ and
$\nabla\rho(z)=-\nu(z)$, Taylor expansion at $z$ gives
\[
 \rho(y)
 =-\nu(z)\cdot(y-z)+O(|y-z|^2)
 =h\bigl(r-w\cdot\nu(z)\bigr)+O_R(h^2)
 =ht+O_R(h^2).
\]
For $a=x,y$, the flow equation gives $\Phi_h(a)-a=hX(z)+O_R(h^2)$.
Indeed, $|a-z|\le C_Rh$ and the trajectory moves by at most $Ch$
for $0\le\tau\le h$, so $X(\Phi_\tau(a))=X(z)+O_R(h)$.
Applying Taylor's formula to $\rho$ at $z$ once more gives
\[
\frac{\rho(\Phi_h(x))}{h}=r-X(z)\cdot\nu(z)+O_R(h)\qquad\text{and}\qquad\frac{\rho(\Phi_h(y))}{h}=t-X(z)\cdot\nu(z)+O_R(h).
\]
These estimates are uniform for fixed $R$. If $r>0$, then $d(x)=hr$, so
$h^{-s}u(x)=\psi(x)r^s\to\psi(z)r^s$. If $r\le0$, the zero
extension gives $u(x)=0$. Thus, $h^{-s}u(x)\to\psi(z)g(r)$, where $g(\tau)=\tau_+^s$.\newline 
The same argument applies to $y$, $\Phi_h(x)$ and $\Phi_h(y)$ using their rescaled signed distances. More precisely, for each of
the four points $a$, if $a\in\Omega$, then
\[
 h^{-s}u(a)-\psi(z)g\!\left(\frac{\rho(a)}{h}\right)
 =(\psi(a)-\psi(z))\left(\frac{\rho(a)}{h}\right)^s.
\]
This tends uniformly to zero because $|a-z|\le C_Rh$,
$\rho(a)/h$ is bounded and $\psi$ is uniformly continuous on
$\overline\Omega$. If $a$ lies on or outside the boundary, both
terms on the left vanish. Since $g$ is uniformly continuous on
bounded intervals, we can replace $\rho(a)/h$ by its limiting
value even when that value is zero. Thus
\[
 \begin{aligned}
 h^{-s}u(x)&\longrightarrow\psi(z)g(r),&
 h^{-s}u(y)&\longrightarrow\psi(z)g(t),\\
 h^{-s}u_h(x)&\longrightarrow\psi(z)g(r-X(z)\cdot\nu(z)),&
 h^{-s}u_h(y)&\longrightarrow\psi(z)g(t-X(z)\cdot\nu(z)),
 \end{aligned}
\]
uniformly on the rescaled region. For points inside $\Omega$ this follows from $u=\psi d^s$ and the
uniform continuity of $\psi$ on $\overline\Omega$. At exterior points
the zero extension of $u$ agrees with the positive-part convention in
$g$. 

In the change $x=z-hr\nu(z)$, differentiation in $r$ gives
$-h\nu(z)$, while differentiation along a tangent vector $v$ gives
$v-hrD\nu(z)v$. The normal factor in the volume element is therefore
$h$ and the tangential Jacobian is
$\bigl|\det_{T_z\partial\Omega}(I-hrD\nu(z))\bigr|=1+O(h|r|)$.
For each fixed $x$, the change $y=x+hw$ gives $dy=h^N\,dw$.
Combining the two changes yields
\[
 dx\,dy=h^{N+1}(1+O(h|r|))\,d\sigma(z)\,dr\,dw.
\]

The kernel contributes a factor $h^{-N-sp}$ and homogeneity gives
$B_p(\delta u,\delta u_h)=h^{sp}B_p(h^{-s}\delta u,h^{-s}\delta u_h)$.
Including the factor $1/h$ and the volume factor $h^{N+1}$, we have
$h^{-1}h^{-N-sp}h^{sp}h^{N+1}=1$. The transformed integral is therefore
\[
\frac1h\iint_{\mathcal C_{h,R,\varepsilon}}\!\!
 B_p(\delta u,\delta u_h)K_0\,dx\,dy=\frac{C_{N,s,p}}2\!
 \int_{\partial\Omega}\int_{-R}^R\int_{\varepsilon<|w|<R}\!\!
 \frac{B_p(h^{-s}\delta u,h^{-s}\delta u_h)}{|w|^{N+sp}}
 (1+O(h|r|))\,dw\,dr\,d\sigma(z).
\]
The differences on the right are evaluated at $x=z-hr\nu(z)$ and
$y=x+hw$. For fixed $R$ and $\varepsilon$, the integrand is uniformly
bounded on the integration region, so dominated convergence gives
\begin{multline}\label{compact-limit}
 \lim_{h\downarrow0}\frac1h
 \iint_{\mathcal C_{h,R,\varepsilon}}
 B_p(\delta u,\delta u_h)K_0\,dx\,dy\\
 =\frac{C_{N,s,p}}2
 \int_{\partial\Omega}|\psi(z)|^p
 \int_{-R}^R\int_{\varepsilon<|w|<R}
 \frac{B_p(g(r)-g(t),\,g(r-X\cdot\nu)-g(t-X\cdot\nu))}
 {|w|^{N+sp}}\,dw\,dr\,d\sigma(z).
\end{multline}
In this expression $X,\nu$ are evaluated at $z$.

\emph{2. Integration parallel to the tangent plane.} Fix $z$ and rotate coordinates so that $\nu(z)=e_N$.
Write $w=(w',w_N)$, where $w'\in\R^{N-1}$ lies in the tangent plane.
Then $t=r-w_N$ and only the kernel depends on $w'$. To integrate it
over the tangent plane, use the gamma integral for $w_N\ne0$.
\[
 (|w'|^2+w_N^2)^{-\frac{N+sp}{2}}
 =\frac1{\Gamma\left(\frac{N+sp}{2}\right)}
 \int_0^\infty
 \tau^{\frac{N+sp}{2}-1}e^{-\tau(|w'|^2+w_N^2)}\,d\tau.
\]
Tonelli and the product Gaussian integral
$\int_{\R^{N-1}}e^{-\tau|w'|^2}\,dw'=(\pi/\tau)^{(N-1)/2}$
therefore give
\begin{equation}\label{tangential}
\begin{split}
 \int_{\R^{N-1}}(|w'|^2+w_N^2)^{-\frac{N+sp}{2}}\,dw'&=\frac{\pi^{\frac{N-1}{2}}}{\Gamma\left(\frac{N+sp}{2}\right)}\int_0^\infty \tau^{\frac{1+sp}{2}-1}e^{-\tau w_N^2}\,d\tau=\frac{\pi^{\frac{N-1}{2}}\Gamma\left(\frac{1+sp}{2}\right)}{\Gamma\left(\frac{N+sp}{2}\right)}|w_N|^{-1-sp}.
\end{split}
\end{equation}
The formula \eqref{tangential} above is applied only after removing the cutoffs in the nonnegative limiting profile integral. At fixed $R$, the sets $\{\varepsilon<|w|<R\}$ increase as $\varepsilon\downarrow0$; afterwards $\{|r|<R,|w|<R\}$ increase to the entire parameter space as $R\to\infty$. Monotone convergence justifies these passages even before finiteness is checked. The finite value follows from Proposition~\ref{prop:constant}.

For each fixed $z$, put $\ell=X(z)\cdot\nu(z)>0$. The numerator depends on $w_N$ only through $t=r-w_N$. Thus, using \eqref{tangential} the change of variable and \eqref{cellscale}, we get
\[
\begin{aligned}
\int_\R\int_{\R^N}
& \frac{B_p(g(r)-g(r-w_N),g(r-\ell)-g(r-w_N-\ell))}{|w|^{N+sp}}\,dw\,dr\\
&\mkern+200mu=\frac{\pi^{\frac{N-1}{2}}\Gamma\left(\frac{1+sp}{2}\right)}{\Gamma\left(\frac{N+sp}{2}\right)}
 \int_\R\int_\R
 \frac{B_p(g(r)-g(t),g(r-\ell)-g(t-\ell))}{|r-t|^{1+sp}}\,dt\,dr\\
&\mkern+200mu=\frac{\pi^{\frac{N-1}{2}}\Gamma\left(\frac{1+sp}{2}\right)}{\Gamma\left(\frac{N+sp}{2}\right)}\,\ell K_{s,p}.
\end{aligned}
\]
Multiplying by
$C_{N,s,p}|\psi(z)|^p/2$ and using \eqref{cdef} gives
$c_{s,p}|\psi(z)|^pX(z)\cdot\nu(z)$. Its integral in $z$ is the
boundary term in \eqref{lower} and \eqref{flowlimit}. For every fixed
$R,\varepsilon$, nonnegativity yields
\[
 \liminf_{h\downarrow0}\frac{R_h}{h}
 \ge\lim_{h\downarrow0}h^{-1}\iint_{\mathcal C_{h,R,\varepsilon}}B_p(\delta u,\delta u_h)K_0.
\]
Letting $\varepsilon\downarrow0$ and then $R\to\infty$ proves \eqref{lower}.
The remaining steps use the equation to exclude additional contributions.

\emph{3. Removing the small-separation cutoff.}
For the upper bound, use now the regularity of a weak solution.
Write $\mathcal C_{h,R,0}$ for \eqref{core} with $0<|x-y|<Rh$.
At fixed $R$, we estimate the omitted contribution divided by $h$,
first taking $\limsup_{h\downarrow0}$ and then letting
$\varepsilon\downarrow0$. Because of \eqref{Bupper}, 
\[
C \iint\limits_{\substack{|\rho(x)|<Rh\\0<|x-y|\le\varepsilon h}}B_p(\delta u,\delta u_h)K_0\,dx\,dy\leq  \iint\limits_{\substack{|\rho(x)|<Rh\\0<|x-y|\le\varepsilon h}} \frac{|u(x)-u(y)|^p}{|x-y|^{N+sp}}\,dx\,dy+\iint\limits_{\substack{|\rho(x)|<Rh\\0<|x-y|\le\varepsilon h}} \frac{|u_h(x)-u_h(y)|^p}{|x-y|^{N+sp}}
\]
Polar coordinates $y=x+\tau\theta$, $\theta\in\mathbb S^{N-1}$, give 
\begin{equation}\label{small-energy}
\begin{aligned}
&\iint_{\substack{|\rho(x)|<Rh\\0<|x-y|\le\varepsilon h}}
 \frac{|u(x)-u(y)|^p}{|x-y|^{N+sp}}\,dy\,dx=\int_0^{\varepsilon h}\tau^{-1-sp}\mkern-30mu
 \int\limits_{\{|\rho(x)|<Rh\}}\mkern-15mu\int_{\mathbb S^{N-1}}
 |u(x)-u(x+\tau\theta)|^p\,d\sigma(\theta)\,dx\,d\tau.
\end{aligned}
\end{equation}
For each fixed $\tau$, let us split the integral in $x$ in three regions. On $|\rho(x)|\le2\tau$, a collar of width $2\tau$, the volume is at most $C\tau$ and $|u(x)-u(x+\tau\theta)|^p\le C\tau^{sp}$ by Lemma~\ref{lem:regularity}. This part is bounded by
\begin{equation}\label{small-collar}
 C\int_0^{\varepsilon h}\tau^{-1-sp}\,\tau\,\tau^{sp}\,d\tau
 =C\int_0^{\varepsilon h}1\,d\tau.
\end{equation}
On $\rho(x)<-2\tau$, both points are exterior and the increment is zero. Finally, on deeper interior points, say, $2\tau<\rho(x)<Rh$, both points are interior and $|x-y|=\tau<d(x)/2$. By \eqref{weightedholder},
\[
 |u(x)-u(x+\tau\theta)|^p\le C\tau^{p\beta}d(x)^{-p(\beta-s)}.
\]
Integrating the last factor in normal coordinates $t=d(x)$ gives
\[
 \int_{\{2\tau<d(x)<Rh\}}d(x)^{-p(\beta-s)}\,dx
 \le C\int_{2\tau}^{Rh}t^{-p(\beta-s)}\,dt.
\]
Consequently this third region contributes at most
\begin{equation}\label{small-interior}
 C\int_0^{\varepsilon h}\tau^{-1-sp}\tau^{p\beta}
 \int_{2\tau}^{Rh}t^{-p(\beta-s)}\,dt\,d\tau.
\end{equation}
Adding \eqref{small-collar} and \eqref{small-interior}, with zero contribution from the exterior region, bounds \eqref{small-energy} by
\[\begin{split}
 C\int_0^{\varepsilon h}
 \left[1+\tau^{p(\beta-s)-1}
 \int_{2\tau}^{Rh}t^{-p(\beta-s)}\,dt\right]d\tau&=Ch\int_0^\varepsilon\left[
 1+a^{p(\beta-s)-1}\int_{2a}^R b^{-p(\beta-s)}\,db
 \right]da\\
&\leq C_Rh\,\varepsilon^{\min\{1,p(\beta-s)\}}
 (1+|\log\varepsilon|),
 \end{split}
\]
for $\varepsilon<R/2$ small enough. For $u_h$, change both variables by $\Phi_h$. The two side Lipschitz bounds \eqref{biLip} control the kernel and volume; the collar becomes an $(R+C)h$
collar and the separation is at most $C\varepsilon h$. More precisely, 
writing $\xi=\Phi_h(x)$ and $\eta=\Phi_h(y)$, the transformed numerator becomes $|u(\xi)-u(\eta)|^p$. The kernel and inverse Jacobians are bounded by a fixed multiple of $|\xi-\eta|^{-N-sp}$. The original restrictions imply $|\rho(\xi)|<(R+C)h$ and $|\xi-\eta|<C\varepsilon h$. Enlarging to this set gives the bound already proved with $R+C$ and $C\varepsilon$ in place of $R,\varepsilon$. These constants are fixed while $\varepsilon\downarrow0$. Therefore,
\begin{equation}\label{diagonal}
 \lim_{\varepsilon\downarrow0}\limsup_{h\downarrow0}\frac1h
 \iint_{\substack{|\rho(x)|<Rh\\|x-y|<\varepsilon h}}
 B_p(\delta u,\delta u_h)K_0\,dx\,dy=0.
\end{equation}
\begin{remark}
Only $\beta>s$ was needed here, so this step remains valid throughout $0<s<1$, independently of the range required for the interior comparison.
\end{remark}

\emph{4. Near pairs and the outer sum.}
We call a pair $(x,y)$ near if $x,y\in\Omega$ and
\[
 |x-y|<\kappa\min\{d(x),d(y)\},\qquad \kappa=1/128.
\]
All other pairs are called separated. Take $R$ large enough that $32/R<c$, where $c$ is the constant
in Lemma~\ref{lem:comparison}. We also take $R>C_0$, with $C_0$
as in Lemma~\ref{lem:increments}.

Fix $t\ge Rh$ and consider the layer $\{x\in\Omega:t\le d(x)<2t\}$.
Choose a maximal collection of centers $x_i$ in this layer with
mutual distances at least $t/256$. By maximality, the balls
$B_{t/256}(x_i)$ cover the layer. The balls $B_{t/512}(x_i)$ are
pairwise disjoint and lie in $\{x\in\Omega:t/2<d(x)<3t\}$.
As for \eqref{layer-Lq}, the enlarged layer has volume
at most $Ct$ so the number of centers is at most $Ct^{1-N}$.

Now let $(x,y)$ be a near pair with $t\le d(x)<2t$ and choose
a center $x_i$ such that $x\in B_{t/256}(x_i)$. Since the pair is
near, $|y-x_i| \le |y-x|+|x-x_i|<2\kappa t+t/256<t/32$. Thus both $x$ and $y$ belong to $B_r(x_i)$ with $r=t/32$.
Since $x_i$ lies in the layer, $32r\le d(x_i)<64r$ and $h/r\le32/R<c$.
We can therefore apply Lemma~\ref{lem:comparison} to $B_r(x_i)$. Using \eqref{local-defect}, we obtain
{\small\begin{equation}\label{near_ine}
\begin{split}
 \iint\limits_{\substack{\text{near pairs}\\t\le d(x)<2t}}
 B_p(\delta u,\delta u_h)K_0\,dx\,dy&\le
 \sum_i\iint\limits_{B_r(x_i)^2}B_p(\delta u,\delta u_h)K_0\,dx\,dy\le Ct^{1-N}r^N\left(\frac{h}{r}\right)^{q\alpha}
 \le Ch^{q\alpha}t^{1-q\alpha}.
 \end{split}
\end{equation}}
The disjoint layers with $t=2^kRh$, $k\ge0$, cover the set
$\{x\in\Omega:d(x)\ge Rh\}$. Thus, \eqref{near_ine} gives
\begin{equation} \label{near-tail}
\iint_{\substack{\text{near pairs}\\d(x)\ge Rh}}
B_p(\delta u,\delta u_h)K_0\,dx\,dy\le Ch^{q\alpha}\sum_{k=0}^{\infty}(2^kRh)^{1-q\alpha}=\frac{ChR^{1-q\alpha}}{1-2^{1-q\alpha}}.
\end{equation}
Where we used $q\alpha>1$. After
division by $h$, the bound is uniform in $h$ and tends to zero
as $R\to\infty$.

\emph{5. Separated pairs and the inner sum.}
For such pairs, \eqref{Bupper} and global $C^s$ regularity give
\[
 B_p(\delta u,\delta u_h)K_0
 \le C|x-y|^{-N-sq}(|D_h(x)|^q+|D_h(y)|^q).
\]
The above bound follows since $|\delta u|+|\delta u_h|\le C|x-y|^s$ from Lemma~\ref{lem:regularity} jointly with the Lipschitz bound for the flow and $|\delta D_h|^q\le2^{q-1}(|D_h(x)|^q+|D_h(y)|^q)$. Let us stress that this bound holds for all pairs; the separated-pair property is used next to bound the integral in $y$. By symmetry, the terms involving $|D_h(x)|^q$ and $|D_h(y)|^q$
have the same integral over the separated pairs with $|x-y|>Rh$.
If $D_h(x)\ne0$, then $x\in\Omega$.
If $y$ is outside, $d(x)\le|x-y|$.
If both are inside and the pair is separated, either
$d(x)\le\kappa^{-1}|x-y|$ or
$d(y)\le\kappa^{-1}|x-y|$; in the latter case use
$d(x)\le d(y)+|x-y|$.
Thus in every case $d(x)\le(1+\kappa^{-1})|x-y|$.
For separated pairs with $|x-y|>Rh$, symmetry and integration in $y$
therefore give the bound
\begin{equation}\label{separated-weight}
\begin{split}
 \iint_{\substack{\text{separated pairs}\\ |x-y|>Rh}}
 B_p(\delta u,\delta u_h)K_0\,dx\,dy\!
 &\le\! C\int_\Omega \!|D_h(x)|^q\!\!
 \int_{|x-y|\ge \max\left\{Rh,\frac{d(x)}{1+\kappa^{-1}}\right\}}\!
 |x-y|^{-N-sq}\,dy\,dx\\
 &\quad= \frac{C}{sq}
 \int_\Omega |D_h(x)|^q
  \max\left\{Rh,\frac{d(x)}{1+\kappa^{-1}}\right\}^{-sq}\,dx\\
 &\quad\le C\int_\Omega |D_h(x)|^q\max\{Rh,d(x)\}^{-sq}\,dx.
\end{split}
\end{equation}
Before enlarging the integration region, the separated-pair set with
$|x-y|>Rh$ is symmetric. Interchanging $x,y$ therefore gives the same
bound for the term involving $D_h(y)$. Since $D_h=0$ outside $\Omega$,
the remaining integral is over $\Omega$.

Since $R>C_0$, we split the integral on the right-hand side of
\eqref{separated-weight} into the regions $d(x)<C_0h$,
$C_0h\le d(x)<Rh$ and $d(x)\ge Rh$. On the first region, the weight is $(Rh)^{-sq}$ and the
thin-layer estimate from Lemma~\ref{lem:increments} gives
\[
 \int\limits_{\{d(x)<C_0h\}}\mkern-12mu
 |D_h(x)|^q\max\{Rh,d(x)\}^{-sq}\,dx
=(Rh)^{-sq\mkern-12mu}\int\limits_{\{d(x)<C_0h\}}\mkern-12mu|D_h(x)|^q\,dx
 \le C(Rh)^{-sq}h^{sq+1}
 =ChR^{-sq}.
\]

On the second region, the weight is still $(Rh)^{-sq}$. We cover
$\{C_0h\le d(x)<Rh\}$ by the layers $t\le d(x)<2t$ with
$t=2^kC_0h$ with $k\in\mathcal{K}=\{k\geq0: 2^kC_0h<Rh\}$. The last layer may extend past $Rh$ and we
include it to obtain an upper bound. By \eqref{layer-Lq} and $1+q(s-\alpha)>0$,
\begin{align*}
 \int\limits_{\{C_0h\le d(x)<Rh\}}
 |D_h(x)|^q\max\{Rh,d(x)\}^{-sq}\,dx
 &\le (Rh)^{-sq}
 \sum_{k\in\mathcal{K}}
 \int\limits_{\{2^kC_0h\le d(x)<2^{k+1}C_0h\}}|D_h(x)|^q\,dx\\
 &\quad\le Ch^{q\alpha}(Rh)^{-sq}
 \sum_{k\in\mathcal{K}}
 (2^kC_0h)^{1+q(s-\alpha)}\\
 &\quad\le Ch^{q\alpha}(Rh)^{-sq}(Rh)^{1+q(s-\alpha)}
 =ChR^{1-q\alpha}.
\end{align*}
On the third region, the weight is $d(x)^{-sq}$. On each layer
$t\le d(x)<2t$ with $t=2^kRh$, it is bounded above by $t^{-sq}$.
Thus \eqref{layer-Lq} yields
\begin{align*}
 \int\limits_{\{d(x)\ge Rh\}}
 |D_h(x)|^q\max\{Rh,d(x)\}^{-sq}\,dx
 &=\sum_{k=0}^{\infty}
 \int_{\{2^kRh\le d(x)<2^{k+1}Rh\}}
 |D_h(x)|^q d(x)^{-sq}\,dx\\
 &\le\sum_{k=0}^{\infty}(2^kRh)^{-sq}
 \int_{\{2^kRh\le d(x)<2^{k+1}Rh\}}|D_h(x)|^q\,dx\\
 &\le Ch^{q\alpha}\sum_{k=0}^{\infty}(2^kRh)^{1-q\alpha}
 \le Ch^{q\alpha}(Rh)^{1-q\alpha}
 =ChR^{1-q\alpha}.
\end{align*}

Adding the three estimates in \eqref{separated-weight}, we obtain
\begin{equation}\label{sep}
 \iint_{\substack{\text{separated pairs}\\ |x-y|>Rh}}
 B_p(\delta u,\delta u_h)K_0\,dx\,dy
 \le Ch\bigl(R^{-sq}+R^{1-q\alpha}\bigr).
\end{equation}

\emph{6. Taking the limits.} If $|x-y|>Rh$, separated pairs were treated in \eqref{sep} above; a near pair with
this separation must have $d(x)>Rh/\kappa$ and is covered by
\eqref{near-tail}.
If $|\rho(x)|>Rh$ but $|x-y|\le Rh$, an exterior $x$ forces $y$
to be exterior, so both functions vanish at both points.
For interior $x$, near pairs are covered by \eqref{near-tail}.
Separated pairs have $|x-y|\ge c\,d(x)>cRh$ and the preceding
separated-pair estimate applies with $R$ replaced by $cR$, for $R$
large enough that $cR>C_0$. Consequently,
\begin{equation}\label{tightness}
 \limsup_{h\downarrow0}\frac1h
 \iint_{\R^{2N}\setminus\mathcal C_{h,R,0}}
 B_p(\delta u,\delta u_h)K_0\,dx\,dy
 \le C(R^{-sq}+R^{1-q\alpha}).
\end{equation}
On the other hand, by \eqref{pair-decomposition}, we have
\begin{align*}
 \frac{R_h}{h}
 -\frac1h\iint_{\mathcal C_{h,R,\varepsilon}}
 B_p(\delta u,\delta u_h)K_0\,dx\,dy=&
 \frac1h\iint_{\substack{|\rho(x)|<Rh\\0<|x-y|\le\varepsilon h}}
 B_p(\delta u,\delta u_h)K_0\,dx\,dy\\
 &+
 \frac1h\iint_{\R^{2N}\setminus\mathcal C_{h,R,0}}
 B_p(\delta u,\delta u_h)K_0\,dx\,dy.
\end{align*}
The first integral in the right hand side is bounded by \eqref{diagonal} while \eqref{tightness}
bounds the second. Taking the upper limit as $h\downarrow0$
and applying \eqref{compact-limit} to the integral over
$\mathcal C_{h,R,\varepsilon}$, we obtain
\begin{align*}
 \limsup_{h\downarrow0}\frac{R_h}{h}
 &\le
 \frac{C_{N,s,p}}2
 \int_{\partial\Omega}|\psi(z)|^p
 \int_{-R}^R\int\limits_{\varepsilon<|w|<R}
\mkern-10mu  \frac{B_p(g(r)-g(t),\,g(r-X\cdot\nu)-g(t-X\cdot\nu))}
 {|w|^{N+sp}}\,dw\,dr\,d\sigma(z)\\
 &\quad+
 C_R\varepsilon^{\min\{1,p(\beta-s)\}}
      (1+|\log\varepsilon|)
 +C\bigl(R^{-sq}+R^{1-q\alpha}\bigr).
\end{align*}

We first let $\varepsilon\downarrow0$ with $R$ fixed. Since
$\beta>s$, the term involving $\varepsilon$ tends to zero.
In the limiting profile integral from \eqref{compact-limit},
the region of integration in $w$ increases from
$\{\varepsilon<|w|<R\}$ to $\{0<|w|<R\}$. The integrand is
nonnegative, so this passage follows by monotone convergence.

We then let $R\to\infty$. Since $sq>0$ and $q\alpha>1$, both
$R^{-sq}$ and $R^{1-q\alpha}$ tend to zero. The regions
$\{|r|<R,\ 0<|w|<R\}$ increase to the full profile region,
up to the null set $w=0$. By monotone convergence and the
calculation in Step~2, we conclude that
\[
 \limsup_{h\downarrow0}\frac{R_h}{h}
 \le c_{s,p}\int_{\partial\Omega}
 |\psi(z)|^p X(z)\cdot\nu(z)\,d\sigma(z).
\]
The lower bound \eqref{lower} gives the reverse inequality
for the lower limit. This proves \eqref{flowlimit}.
\end{proof}

\begin{proof}[Proof of Theorem~\ref{thm:general}]
Assume first that $X$ is outward. The exact cancellation
\eqref{flow-exact} was obtained by testing with $D_h$.
Because $f$ is Lipschitz on the closed interval containing the values
of the zero extension of $u$,
\[
 |F(b)-F(a)-f(a)(b-a)|\le C|b-a|^2.
\]
Using boundedness of $D_h$, $q\le2$ and \eqref{global-Lq}, we obtain
\[
 |T_h|\le C\|D_h\|_2^2\le C\|D_h\|_q^q
 \le Ch^{q\alpha}=o(h).
\]

We now differentiate the right side of \eqref{flow-exact}.
For the energy use $x=\Phi_{-h}(\xi)$,
$y=\Phi_{-h}(\eta)$, then rename the dummy integration variables:
\[
 \frac{C_{N,s,p}}2[u_h]_{s,p}^p
 =\iint_{\R^{2N}}|u(x)-u(y)|^pK_{-h}(x,y)\,dx\,dy.
\]
Since $F(0)=0$ and $\operatorname{supp}u_h\subset\Omega$ for small
$h>0$, we have $F(u_h)=0$ outside $\Omega$. We may therefore extend
the integral to $\R^N$ and use $u_h=u\circ\Phi_h$. Thus,
\[
 \int_\Omega F(u_h(x))\,dx=\int_{\R^N}F(u(\Phi_h(x)))\,dx=\int_{\R^N}F(u(y))j_{-h}(y)\,dy=\int_\Omega F(u(y))j_{-h}(y)\,dy.
\]
Subtracting the corresponding identities at $h=0$, the exact energy identity becomes
\[
\begin{aligned}
 \frac{R_h-T_h}{h}
 &=\frac1p\iint_{\R^{2N}}|\delta u|^p
       \frac{K_{-h}-K_0}{h}\,dx\,dy
 -\int_\Omega F(u)\frac{j_{-h}-1}{h}\,dx.
\end{aligned}
\]
The inverse flow has
$\partial_hK_{-h}|_{h=0}=-K_0L_X$ and
$\partial_hj_{-h}|_{h=0}=-\operatorname{div}X$.
The uniform bound in \eqref{kernelbounds} and finiteness of the energy
justify the use of the dominated convergence theorem. Indeed, the two absolute integrands are bounded respectively by $C|\delta u|^pK_0$, integrable because $u\in W^{s,p}(\R^N)$ and $C|F(u)|\mathbf1_\Omega$, integrable because $u$ is bounded and $\Omega$ has finite volume. Passing to the limit gives
\begin{equation}\label{first-variation}
 \lim_{h\downarrow0}\frac{R_h-T_h}{h}
 =\int_\Omega F(u)\operatorname{div}X\,dx
 -\frac1p\iint_{\R^{2N}}|\delta u|^pK_0L_X\,dx\,dy.
\end{equation}
Combine $T_h=o(h)$ with Proposition~\ref{prop:boundary}.
Since $K_0=C_{N,s,p}/(2|x-y|^{N+sp})$, we conclude
\eqref{general-identity} for outward fields. To extend the identity to a general field $X$, let $U$ be a
tubular neighborhood of $\partial\Omega$ and let $\pi(x)$ denote
the nearest boundary point to $x\in U$. Choose
$\chi\in C_c^\infty(U)$ with $\chi=1$ near $\partial\Omega$ and set $V(x)=2\chi(x)\nu(\pi(x))$ if $x\in U$ and $V(x)=0$ otherwise. Then $V\in C_c(\R^N;\R^N)$ and $V(z)=2\nu(z)$ for every
$z\in\partial\Omega$. Let $\eta_\varepsilon$ be a nonnegative smooth mollifier of mass
one supported in $B_\varepsilon$. By uniform continuity of $V$,
choose $\varepsilon$ so that $Z=\eta_\varepsilon*V$ satisfies
$\|Z-V\|_\infty<1$. Then $Z\in C_c^\infty(\R^N;\R^N)$ and
$Z(z)\cdot\nu(z)\ge V(z)\cdot\nu(z)-\|Z-V\|_\infty>1$ for every $z\in\partial\Omega$.

For a general vector field $X$ now take $M>\|X\cdot\nu\|_{L^\infty(\partial\Omega)}$. On
$\partial\Omega$ we have
\[
 (MZ)\cdot\nu>0,
 \qquad
 (X+MZ)\cdot\nu
 \ge -\|X\cdot\nu\|_{L^\infty(\partial\Omega)}+M>0.
\]
Thus, we can apply the identity to $X+MZ$ and to $MZ$.
Since both sides of the identity are linear in the vector field,
subtracting the two equalities gives the identity for
$X=(X+MZ)-MZ$. Taking $X(x)=x$ proves
\eqref{pohozaev}.
\end{proof}

\begin{samepage}
\begin{corollary}\label{cor:nonexistence}
Assume \eqref{new-range}, $sp<N$ and $x\cdot\nu>0$ on
$\partial\Omega$. If $m\ge\max\{2,Np/(N-sp)\}$, there is no bounded
weak solution of
\[
 (-\Delta)^s_pu=|u|^{m-2}u\quad\text{in }\Omega,\qquad
 u=0\quad\text{in }\Omega^c
\]
whose weighted trace is nonzero on a set of positive surface measure.
\end{corollary}
\end{samepage}
\begin{proof}
The condition $m\ge2$ makes the nonlinerity locally Lipschitz.
Since $F(u)=|u|^m/m$, the left side of \eqref{pohozaev} is nonpositive while the right side is strictly positive by the trace assumption,
$x\cdot\nu>0$ and $c_{s,p}>0$, a contradiction.
\end{proof}

\section*{Acknowledgment}
\noindent A. Ortega is partially funded by the Vicerrectorado de Investigación, Transferencia y Divulgación Científica of the Universidad Nacional de Educación a Distancia under the Plan de Promoción de Investigación through the research project \textit{Proyectos Investigación Tipo A. EDP no locales con dato de frontera mixto Dirichlet-Neumann}, Ref. 2026/00131/001.

\end{document}